\documentclass[1p]{elsarticle}
\usepackage{lineno,hyperref}
\usepackage{amsmath}
\usepackage{amsfonts}
\usepackage{amsthm}
\usepackage{upgreek}
\usepackage{graphicx}
\usepackage{subfigure}
\usepackage{appendix}
\usepackage{enumitem}
\journal{arXiv }
\modulolinenumbers[5]

\usepackage[font=small,labelfont=bf,labelsep=space]{caption}
\def\R{\mathbb R} 
\def\T{\mathcal{T}}
\def\D{\mathcal{D}}
\def\E{\mathcal{E}}
\def\d{\text{\rm d}}

\def\m{\text{\rm m}}
\numberwithin{theo}{section} 

\numberwithin{lemm}{section} 

\numberwithin{rk}{section} 
\newtheorem{Def}{Definition}[section]
\numberwithin{Def}{section} 

\numberwithin{prop}{section} 

\begin{document}

\begin{frontmatter}

\title{Linearized Implicit Methods Based on a Single-Layer
Neural Network: Application to Keller-Segel Models }


\author{M. Benzakour Amine 
\corref{mycorrespondingauthor}}
\ead{benzakouramine.m@ucd.ac.ma
}
\cortext[mycorrespondingauthor]{Corresponding author}


\address{     D\'epartement de Math\'ematiques, Facult\'e des Sciences d'El Jadida, Universit\'e Chouaib Doukkali, 24000 El Jadida, Morocco  }

\numberwithin{equation}{section}
\begin{abstract}
This paper is concerned with numerical approximation of some two-dimensional Keller-Segel chemotaxis models, especially those generating pattern formations.   The numerical resolution of such nonlinear parabolic-parabolic or parabolic-elliptic systems of partial differential equations consumes a significant computational time when solved with fully implicit schemes. Standard linearized semi-implicit schemes,  however, require reasonable computational time, but suffer from lack of accuracy. In this work, two methods based  on a single-layer neural network are developed to build linearized implicit schemes: a basic one called the each step training linearized implicit (ESTLI) method  and a more efficient one, the selected steps training linearized implicit (SSTLI) method. The proposed schemes make use also of  a spatial finite volume method with a hybrid difference scheme approximation for convection-diffusion fluxes. Several numerical tests are performed to illustrate the accuracy, efficiency and robustness of the proposed methods. Generalization of the developed methods  to   other nonlinear partial differential equations  is straightforward.

\end{abstract}
\end{frontmatter}

\section{Introduction}











One of the fundamental  characteristics  of living organisms is their ability to move in response to external signals. Chemotaxis, which is the oriented motion induced by  chemical gradients,  has a critical role in many biological and medical fields. It is encountered  in several  self-organization phenomena such as: aggregation of bacteria, embryogenesis, ecology,  tumor growth, etc.

In 1970, Keller and Segel \cite{Keller1} introduced a chemotaxis model to describe the aggregation of slim molds. Since then, several variants of the  Keller-Segel model have been   developed  to model  phenomena in a variety of fields where chemotaxis is involved (see for more details \cite{Painter}). The model is still very popular, it reads in a generalized form

 \begin{equation}
\left\{  
 \begin{aligned}
 &\partial_t u = \nabla \cdot \left( D_u \nabla u -  u\chi(u,c)\, \nabla c \right)+f(u) ,
 \\
 &\partial_t c = D_c \Delta c +g(u,c),
   \end{aligned}
 \right.
 \label{(1)}
\end{equation}

where $u(x,t)$ denotes the cell density at spatial position $x$ and time $t$, $c(x,t)$ is the   chemoattractant concentration, $D_u$ and $D_c$ respectively define cells diffusion and chemical diffusion coefficients,  $\chi(u,c)$ is the chemotactic sensitivity, $f(u)$ is a function describing cell growth and death, and $g(u,c)$ is a kinetic function that describes production and degradation of the chemical signal. 

Several numerical methods  for Keller-Segel systems have been proposed and analyzed. The interested reader may be referred  for instance to the introduction of \cite{Akhmouch1}, and to \cite{Jungel,Huang,Xiao,Zhang1,Sulman,Zhang2,Liu,Li,Dehghan} for more recent works.

In this paper, we focus on three models derived from \eqref{(1)}. The first one has been proposed in the framework of pattern formation in embryology, the model reads:
\begin{equation}
\left\{
\begin{aligned}
 &\partial_t u = D_u \Delta u - \chi \nabla \cdot (  \, u \nabla c)\quad \text{in}\ \Omega_T ,
 \\
&0= \Delta c +\frac{u}{u+1} - c \quad \text{in}\ \Omega_T,
 \end{aligned}
 \right.
 \label{(2)}
\end{equation}
with  zero-flux boundary and initial conditions
\begin{equation}
 \nabla u \cdot \nu=\nabla c \cdot \nu=0 \quad
\text{on} \ \partial \Omega \times (0,T),\quad u(.,0)=u_0   \quad \text{in} \ \Omega .
\label{(3)}
\end{equation}
Here and in the rest of this paper, $D_u$ and $\chi$ are positive constants, $\Omega_T:=\Omega \times (0,T)$,  $\Omega \subset\R^2$ is an open bounded   subset, $T>0$ denotes the final time, $\nu$ is the outward unit-normal on the boundary $\partial \Omega$, and  $u_0$ is a nonnegative function.

The above  model is a simpler version of the original parabolic-parabolic model proposed by Oster and Murray  \cite{Oster1}, in which the second equation of the system is elliptic, using the reasonable assumption that  the chemoattractant diffuses much faster than  cells. In \cite{Murray1}, it has been demonstrated that the model is able to generate stripe patterns similar to that  of embryonic American alligators.

In \cite{Akhmouch2}, a semi-implicit finite volume scheme with an additional correction term has been proposed  by Akhmouch and the author to solve the system \eqref{(2)}--\eqref{(3)}. Moreover, a convergence analysis has been performed in the case of a nonnegative initial cell density $u_0 \in L^2(\Omega)$, and it has been established that the numerical solution of the proposed scheme converges to a weak solution of \eqref{(2)}--\eqref{(3)} under some conditions. In this paper, a particular attention will be given to this analysis.

 
The second Keller-Segel system considered in this paper is the following chemotaxis-growth model:

\begin{equation}
\left\{
\begin{aligned}
 &\partial_t u = D_u \Delta u - \chi \nabla \cdot (  \, u \nabla c)+ f(u)\quad \text{in}\ \Omega_T ,
 \\
&\partial_t c= \Delta c -\lambda c +u  \quad \text{in}\ \Omega_T,
 \end{aligned}
 \right.
 \label{(4)}
\end{equation}
with zero-flux boundary conditions and suitable initial conditions. Herein, $\lambda$ is a positive constant and $f(u)$ is a polynomial growth term. This model has the ability to reproduce different  bacterial pattern formations through the coupled effect of chemotaxis and growth (see, e.g., \cite{Aida,Strehl1,Strehl2,Akhmouch1,Akhmouch3}). It is well known that when $f(u)=0$, the solution of the system (called minimal model) may blow-up in finite time for  initial cell density  with sufficiently large mass, which is an interesting mathematical phenomenon. Several works were devoted to study this system (see for more details \cite{Horstmann1,Hillen1,Perthame1}). However, from biological point of view,  the occurrence of finite-time blow-up is considered as a pathological behavior of the model.  An adequate choice of the reaction term $f(u)$ can prevent this blow-up \cite{Winkler}.

The last Keller-Segel system considered in the current paper reads:

\begin{equation}
\left\{
\begin{aligned}
 &\partial_t u = D_u \Delta u -  \nabla \cdot (  \, u \chi(u) \nabla c)\quad \text{in}\ \Omega_T ,
 \\
&\partial_t c= \Delta c - c +u  \quad \text{in}\ \Omega_T,
 \end{aligned}
 \right.
 \label{(5)}
\end{equation}
endowed  with zero-flux boundary conditions and  initial conditions. In the above model, $\chi$ is a function which satisfies the conditions: $\chi(0)>0$,  there exists $\bar u>0$ such that $\chi(\bar u)=0$ and 
$\chi(u)>0$ for $0<u<\bar u$. The system \eqref{(5)},  called volume-filling chemotaxis model, has been introduced  by Hillen and Painter in \cite{Hillen2}, where the existence of solutions that are bounded is proved. The adopted choice of chemotactic sensitivity implies that cells stop to accumulate when the threshold $\bar u$ is reached by $u$, which prevents overcrowding. Numerical simulations of the model are performed in  \cite{Hillen3}, and in \cite{Andreianov1,Chamoun1} for more general diffusion term.

The three chemotaxis models presented are  highly nonlinear coupled systems due to the chemotaxis term $\chi \nabla \cdot (  \, u \nabla c)$ ( $\nabla \cdot (  \, u \chi(u) \nabla c)$ for the volume-filling model). Moreover, the cell density reaction term of the second equation of \eqref{(2)}, the growth term of the first equation of \eqref{(4)}, and the convective term of the first equation of \eqref{(5)}, are also nonlinear. These nonlinearities are an important factor for the choice of the time discretization approach when numerically solving such systems. The most used approaches in the literature are: explicit, implicit and semi-implicit schemes. The main drawback of explicit schemes is the too restrictive CFL condition necessary to ensure  stability of the solution. For implicit schemes, the CFL constraint is avoided. But, on the other hand, a large nonlinear algebra system  must be solved at each time-step, which is too time-consuming. The third time discretization which is also widely used for solving nonlinear systems is the semi-implicit discretization. Semi-implicit schemes are linearized implicit schemes, so they only require to solve decoupled linear algebra systems at each time level. Regarding numerical methods for Keller-Segel systems, Euler semi-implicit scheme is used in several works to obtain fully linear implicit  schemes (see, e.g.,  \cite{Strehl2,Akhmouch3,Liu,Saito1,Saito2}). In other works, it has only been  used to avoid the severe nonlinearity caused by the chemotaxis term and thus obtain decoupled nonlinear schemes \cite{Andreianov1,Chamoun1,Ibrahim1}. As proved in the convergence analysis carried out by Zhou and Saito \cite{Saito2}   for the parabolic-elliptic minimal model, by Akhmouch and the author \cite{Akhmouch3} for  parabolic-elliptic chemotaxis growth model and by Andreianov et al.  for  degenerate  volume-filling chemotaxis model \cite{Andreianov1}, no CFL condition is required to obtain the convergence of  Euler semi-implicit finite volume schemes for such systems. Nonetheless, despite the advantages of such strategy, these schemes generally suffer from a lack of accuracy.

In this work, a novel approach to develop linear implicit schemes is presented. More precisely, two methods relying on a basic single-layer neural network to linearize schemes are proposed: the each step training linearized implicit (ESTLI) method,  and  the selected steps training linearized implicit (SSTLI) method which is based also on a  selected steps training (STP) algorithm. The developed methods  will be used to numerically solve the systems \eqref{(2)}, \eqref{(4)} and \eqref{(5)}. The objective is to build a method which has the following attractive benefits:

\begin{itemize}
\item[$\bullet$]   It must only require  to solve decoupled linear systems at each time-step.
\item[$\bullet$]  It  should be significantly more accurate than semi-implicit method.
\item[$\bullet$] In terms of computational cost, the method should be nearly as efficient as semi-implicit method, specially for time-consuming numerical experiments.
\item[$\bullet$] As semi-implicit method, it is required to be stable even for large time-step sizes.   
\item[$\bullet$] It should be  easily applied  to any equation or system which can be solved by semi-implicit scheme.
\end{itemize}

Concerning the spatial approximation, a finite volume scheme similar to that of Akhmouch and Benzakour  \cite{Akhmouch2} is adopted. We mention that finite volume method has been first introduced by Filbet \cite{Filbet1} for Keller-Segel models, where a fully implicit finite volume scheme has been proposed for the parabolic-elliptic minimal model.



As mentioned above, the approach presented in this article uses artificial neural networks (ANN) only as a linearization technique in a numerical scheme. We mention, however, that several interesting works were devoted to use ANN as a main tool to approximate solutions of partial differential equations, so we briefly review some of them. In \cite{Lagaris1}, multi-layer perceptron was used by Lagaris et al. to solve boundary value problems.   Lagaris et al.   \cite{Lagaris2} developed also a method based on multi-layer perceptron and radial basis function networks. In \cite{Shirvany}, a novel method was presented using the same networks, with application to the nonlinear Schrodinger  equation. We note that several works proposed radial basis function networks    to solve elliptic equations, see for instance \cite{Mai-Duy,Aminataei,Jianyu}. The reliability of other methods based on single-layer Bernstein neural network \cite{Sun} and single-layer Chebyshev neural network \cite{Mall}, was also demonstrated by application to elliptic equations. In \cite{Rudd}, a method combining  Galerkin methods and ANN was developed  for the approximation of the solution of hyperbolic and parabolic partial differential equations. ANN were also used to solve systems of partial differential equations \cite{Beidokhti}, and high dimensional partial differential equations \cite{Sirignano,Weinan}.

 This paper is organized as follows. In Sect. 2, a semi-implicit Euler  scheme in time and finite volume in space for the the initial-boundary value problem \eqref{(2)}--\eqref{(3)} is presented. The convergence of the scheme is also discussed  in this section on the basis of the analysis performed in \cite{Akhmouch2}. Sect. 3 introduces  ESTLI and SSTLI methods to solve \eqref{(2)}--\eqref{(3)}. The STP algorithm is also presented in this section.  Numerical tests applied to the models presented in this section are performed in Sect. 4, in order to compare the efficiency and accuracy of the elaborate schemes with that of semi-implicit schemes.  Finally, Sect. 5 presents a conclusion.

 \section{The semi-implicit finite volume scheme  } 
\label{s2}
\subsection{Definitions and notations}
  \label{s2.1}
 We assume that $\Omega$  is an open bounded  polygonal  subset. An  admissible finite volume mesh of $\Omega$ in the sens of Definition 9.1 in \cite{Eymard1} is given by:

\begin{itemize}
\item[$\bullet$] A family  of control volumes (disjoint open and convex polygons)  denoted by $\T$. 
\item[$\bullet$] A family $\E$ of  edges, where $\E_K$ is the set of edges of the control volume $K\in \T$.
\item[$\bullet$]  A family of points $(x_K)_{K\in\T}$ such that $x_K \in \overline{K}$.  The particularity of an admissible mesh is that the  straight line going through $x_K$ and $x_L$ must be orthogonal to the common  edge of $K$ and $L$ denoted $K|L$.
\end{itemize}

 We denote by m, the  Lebesgue measure in $\R^2$ or $\R$. Denoting by d  the Euclidean distance, and for all $\sigma \in \E_K$,  we define $\tau_\sigma$  by:

 $$
\tau_\sigma =\left\{\begin{array}{ll}
 \dfrac{\m(\sigma)}{\d(x_K,x_L))},  &\quad\mbox{if } \sigma=K|L, \\
  \dfrac{\m(\sigma)}{\d(x_K,\sigma)}, &\quad\mbox{ if }\sigma\subset   \partial \Omega.
  \end{array}\right.
$$

 Let  $N$  be the number of time-steps to reach the final time $T$. We set $t_n =n\Delta t$, where $\Delta t$ is the time-step size: $\Delta t=\dfrac{T}{N}$.

Let $X(\T)$ be the set of functions from  $\Omega$ to $\R$ which are constant over each control volume of the mesh. Let  $1\leq p<\infty$, for $v \in X(\T)$, the classical discrete $L^p$ norm  reads
$$\|\,v\,\|_p = \left(\sum_{K\in\T}\m(K)\,|v_K|^p\right)^{1/p}, $$
where $v(x)=v_K$ for all $x \in K$ and for all $K \in \T$. 
 We define also  the discrete $H^1$ seminorm and the discrete $H^1$ norm: 
 $$|\,v\,|_{1,\T} = \left(\sum_{\sigma\in\E}\tau_\sigma\, 
  |D_\sigma v|^2\right)^{1/2}, 
$$

$$
\|\,v\,\|_{1,\T} = \|\,v\,\|_2 + |\,v\,|_{1,\T} \, ,
$$ 
where for all $\sigma \in \E$, $D_\sigma v=0$ if $\sigma \subset \partial \Omega$ and $D_\sigma v=|v_K-v_L|$ otherwise, with $\sigma=K|L$.

Finally, we define a weak solution of the system  \eqref{(2)}--\eqref{(3)}:
\begin{Def}
A weak solution of   \eqref{(2)}--\eqref{(3)} is a pair of functions $(u,c)\in L^2(0,T;H^1(\Omega))^2 $ which satisfy the following identities for all test functions $\phi\in \D(\Omega_T)$ : 
\begin{gather}
  \int_0^T\int_\Omega\left(u\,\partial_t \phi - D_u \nabla u\cdot\nabla \phi + \chi\, u\nabla c\cdot\nabla \phi \right)\,dxdt
  + \int_\Omega u_0\,\phi(x,0)\,dx = 0,\label{(f1)}  \\
  \int_0^T\int_\Omega\nabla c\cdot\nabla\phi \ dxdt = \int_0^T\int_\Omega \left(\frac{u}{u+1}- c\right)\phi \,dxdt.\,\label{(f2)}
\end{gather}
\end{Def}

 \subsection{Presentation of the scheme }
  \label{s2.2}

 A semi-implicit  finite volume scheme associated to the      problem      \eqref{(2)}--\eqref{(3)} is given by:
\\
\\
for all $K\in \T$ and $n=0,...,N-1$,

\begin{align}
& \m(K)\frac{u^{n+1}_K-u^n_K}{\Delta t}
  - D_u \sum_{\sigma\in\E_K}\tau_\sigma Du_{K,\sigma}^{n+1} \notag
  \\
  &+\chi \sum_{\substack{\sigma\in\E_{K}\\ \sigma=K|L}}\tau_\sigma\left(S\left( Dc_{K,\sigma}^{n+1}\right)u_{K}^{n+1}-S\left( -Dc_{K,\sigma}^{n+1}\right)u_{L}^{n+1}\right)=0, \label{(6)}
 \\
  & -\sum_{\sigma\in\E_K}\tau_\sigma\, Dc^{n+1}_{K,\sigma}
  +\m(K)\left( c_K^{n+1}-  \frac{u_K^{n}}{u_K^{n}+1}\right)=0, \label{(7)}
\end{align}
endowed with the discrete initial condition
\begin{equation}
 u^0_K = \frac{1}{\m(K)}\int_K u_0(x)\,dx.\label{(8)}
 \end{equation}
  The function $S$ is defined by
  \begin{equation}
 S(x) = \left\{\begin{array}{ll}
0, &\quad\mbox{if }x<2 \left(\varepsilon-D_u \right)/\chi,\\ 
x, &\quad\mbox{if } x>2 \left( D_u-\varepsilon \right)/\chi,\\
  \dfrac{x}{2}, &\quad\mbox{otherwise},
  
   \end{array}\right. \\ \label{(S)}
\end{equation}
where $\varepsilon$ is a small nonnegative constant ($\varepsilon <<D_u$), and $Dv_{K,\sigma}^n =v_L^n-v_K^n$ if $\sigma=K|L \not \subset \partial \Omega$ and $0$ otherwise.

 In the above scheme,  $u^n_K$ and $c^n_K$ denote respectively the approximations of the mean value of $u$ and $c$ on $K$ at time $t_n$. As we can see, the difference between this scheme and a fully implicit one is the discretization of the reaction term $\frac{u}{u+1}$, which is evaluated at the previous time-step.   This strategy allows us to build a  decoupled linear scheme:  at each time-step, we solve \eqref{(7)} to obtain $c^{n+1}_K$ and then, we solve \eqref{(6)} to compute $u^{n+1}_K$ . The  discretization used for the chemotaxis term is identical to  the first order upwind scheme when  $\vert Dc_{K,\sigma}^{n+1}\vert> 2 \left(D_u-\varepsilon\right)/\chi$, and to the second order central difference scheme  otherwise. It is equivalent to that of Spalding \cite{Spalding1}  when $\varepsilon=0$.

 \subsection{Convergence of the scheme }
  \label{s2.3}

The scheme \eqref{(6)}--\eqref{(8)} is  similar to that of \cite{Akhmouch2}, where a corrected decoupled scheme is proposed for the problem \eqref{(2)}--\eqref{(3)}. The only difference between the two schemes is a bounded term which was added to the second equation of the scheme:

\begin{equation*}
 -\sum_{\sigma\in\E_K}\tau_\sigma\, Dc^{n+1}_{K,\sigma}
  +\m(K)\, c_K^{n+1}= \m(K)\, \frac{u_K^{n}}{u_K^{n}+1}+\beta_n T_K^{n},\label{(c)}
 \end{equation*}
we refer to Section 2.2 of  \cite{Akhmouch2} for the definition of $\beta_n T_K^{n}$. The numerical analysis of the two schemes is quite similar, so based on the results obtained for the corrected decoupled scheme, we briefly discuss the convergence of \eqref{(6)}--\eqref{(8)}.

The existence, uniqueness and nonnegativity of  the finite volume solution  $\{(u_{K}^{n+1},c_{K}^{n+1}),\,K \in \T,\, n=0,...,N-1\} $ of  \eqref{(6)}--\eqref{(8)} can be proved by an M-matrix analysis, following exactly the proof of Proposition 3.1 in \cite{Akhmouch2}.

Now, to prove the convergence of the scheme to a weak solution of  \eqref{(2)}--\eqref{(3)} in the sense of Definition 2.1, we first need a priori estimates on the discrete solutions. For the finite volume  approximation of $c$, $\displaystyle L^{\infty}(\Omega)$ and discrete $\displaystyle L^{\infty}(0,T;H^1(\Omega))$ estimates can be obtained directly from the equation  \eqref{(6)}. We refer to Proposition 4.1 of \cite{Akhmouch1}, where the same equation is analyzed. 

For the cell density finite volume solution, only  $L^2(\Omega)$ estimate  is obtained for the corrected decoupled scheme without establishing the boundedness of the solution, since it is supposed that $u_0 \in L^2(\Omega)$. Two approaches have been proposed to obtain this estimate. The first one, detailed in the proof of  Lemma 4.1 in  \cite{Akhmouch2}, needs that $\varepsilon>0$, and holds also for the semi-implicit scheme \eqref{(6)}--\eqref{(8)}. The main tool of this approach is the discrete  Gagliardo-Nirenberg-Sobolev inequality \cite{Chatard2}, which requires  the  following constraint on the mesh: there exists $\xi>0$ such that
$\d(x_K,\sigma)\geq \xi \,\d(x_K,x_L).$ The second approach (Proposition 4.3 in  \cite{Akhmouch2}) is based on the discrete Gronwall inequality, and requires a time-step condition which can be relaxed in the case of the current semi-implicit scheme.

The $\displaystyle L^{2}(0,T;H^1(\Omega))$ estimate for the cell density can be  obtained by following the proof of Proposition 4.2 in  \cite{Akhmouch2}. The rest of the convergence analysis of the scheme \eqref{(6)}--\eqref{(8)},  including the compactness of the sequence of approximate solutions and the  pass to the limit in the scheme,   can be performed exactly as in Section 5 in \cite{Akhmouch2} (see also \cite{Hillairet1,Filbet1}).

\section{Linearization via a single-layer neural network  }

In this section, for all family of elements $\left( v_k \right)_{K \in \T}$, we denote by $v$ the vector    whose components are the elements of the family.

\subsection{The ESTLI  method}
The idea of the ESTLI method is to replace   \eqref{(7)} in the semi-implicit scheme by the equation:
\begin{equation}
 -\sum_{\sigma\in\E_K}\tau_\sigma\, Dc^{n+1}_{K,\sigma}
  +\m(K)\left( c_K^{n+1}-  \frac{\widetilde u_K^{n+1}}{\widetilde u_K^{n+1}+1}\right)=0,\label{(EST)}
 \end{equation}
 where $\widetilde u_K^{n+1}$ is an approximation of $u_K^{n+1}$, obtained through the use of an ANN. The training data are obtained from the previous values of the finite volume solution, which implies that the numerical solution must be computed by an other method for the first time-steps. 
 
 It is easy to see that the ESTLI scheme consumes more time than the semi-implicit scheme  \eqref{(6)}--\eqref{(7)}. Indeed, at each time step, in addition  to the time required  for forming and solving linear systems, a computational time is needed to train the neural network and to compute $\widetilde u_K^{n+1}$ for all $K\in \T$. As training single-layer neural networks already requires much computational time, and keeping in mind that the scheme must be efficient (see the objectives listed in Sect. 1), hidden layers are avoided. The network proposed is a single-layer feedforward neural network which consists of two inputs, weights, a bias and an output. The activation function will be the identity activation function.
 
\begin{figure}[h!]
\begin{center}
\includegraphics[width=0.85\textwidth,height=0.2\textheight]{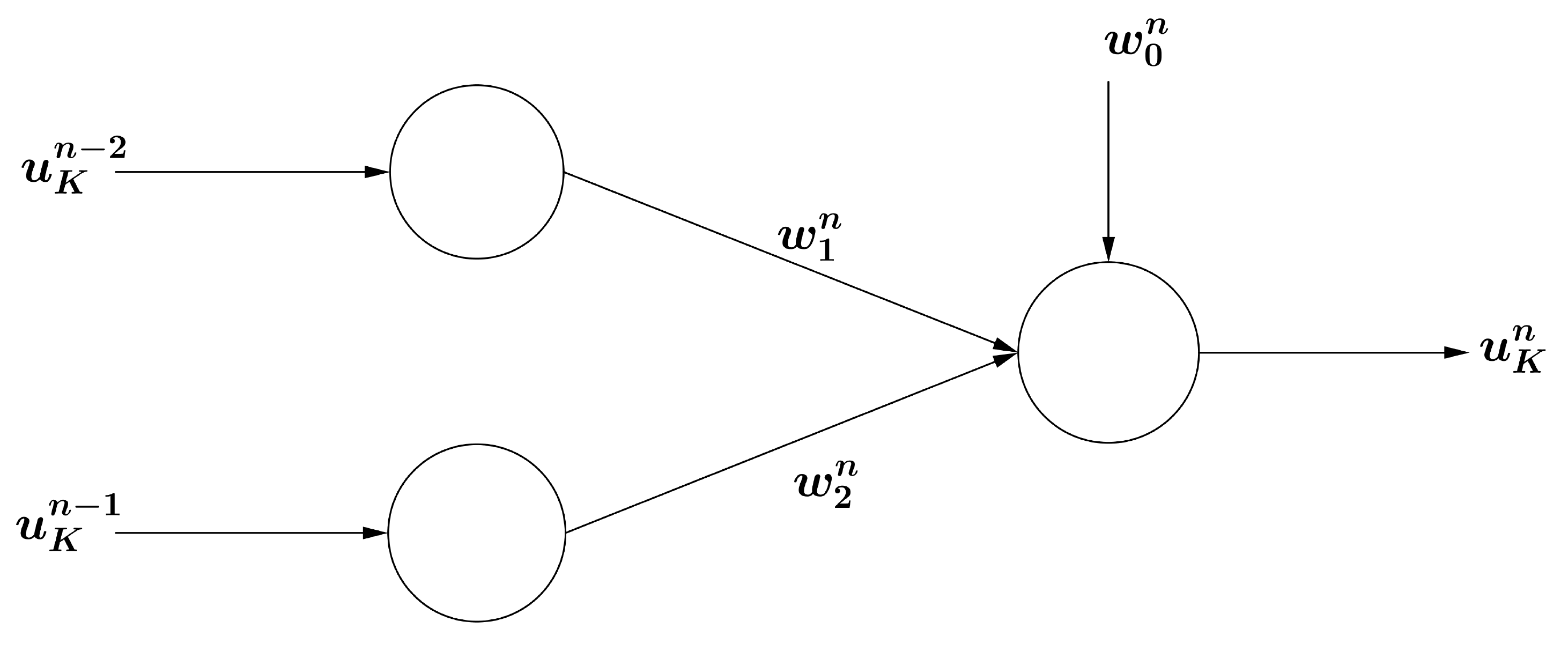}
\caption{Structure of the  neural network at time-step $n+1$, $n \geq 2$ (Training phase) }
\label{Fig.1}
\end{center}
\end{figure}

\begin{figure}[h!]
\begin{center}
\includegraphics[width=0.85\textwidth,height=0.2\textheight]{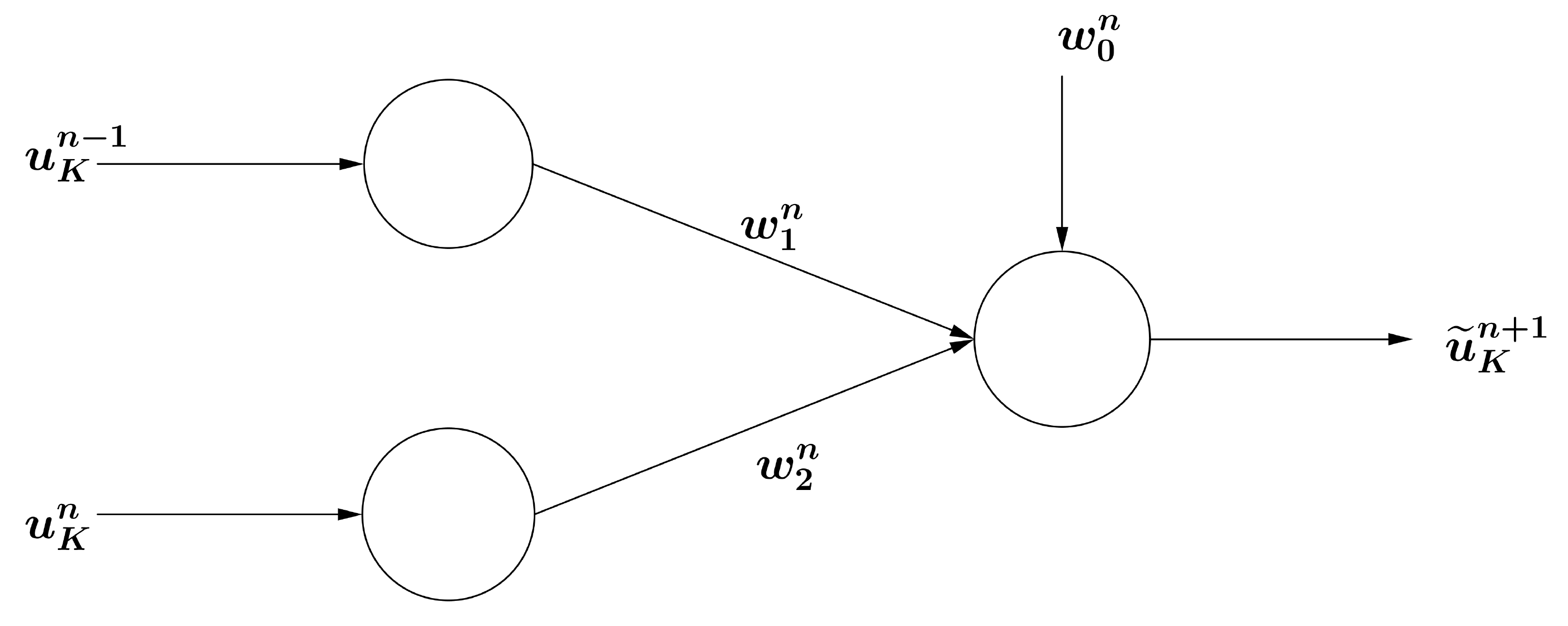}
\caption{Structure of the  neural network at time-step $n+1$, $n \geq 2$ (Prediction phase) }
\label{Fig.2}
\end{center}
\end{figure}

 The structures of the ANN used in the training phase and prediction phase are presented in Figs. \ref{Fig.1} and \ref{Fig.2} respectively. As we can see, at each time-step $n+1$ with $n \geq 2$, the network is trained using the two inputs $u_K^{n-2}$ and $u_K^{n-1}$, and the target output $u_K^n$. The semi-implicit scheme  \eqref{(6)}--\eqref{(7)} will be used for the two first time-steps ($n=0$ and $n=1$) to obtain the numerical solution. The values of the discrete solution at each control volume $K \in \T$ are used in the training process. After completing the training, the final weights $(w_1^n,w_2^n)$ and the bias $w_0^{n}$  are used to compute $\widetilde u^{n+1}$, so for all $K \in \T$:   
 \begin{equation}
 \widetilde u_K^{n+1}=w_1^{n} u_K^{n-1}+w_2^{n} u_K^{n}+w_0^{n}. \label{(20)}
 \end{equation}
 
 The  training process adopted  is a standard one: the weights and the bias are updated until the error reaches the desired  threshold $\epsilon$. Denoting by $w^{n,k}$ the vector of  weights and the bias  at the time step  $n+1$ and the iteration $k$ (related to the updating process), the update rule has the following form:
 \begin{equation}
 w^{n,k+1}=w^{n,k}+  \Updelta w^{n,k}, \label{(21)}
 \end{equation}
where the expression of $\Updelta w^{n,k}$ is determined according to the method selected as learning algorithm:  gradient descent method, conjugate gradient method, Levenberg-Marquardt method, etc. We refer to \cite{Yu,Tan} for an overview of update rules related to several optimization methods. The error function adopted to evaluate the training process is the mean squared error function $E^{n,k}$:
\begin{equation}
E^{n,k}=\dfrac{1}{\mbox{card}(\T)}\sum_{K\in\T} \left(u^{n}_K-\hat{u}^{n,k}_K   \right)^2, \label{(22)}
\end{equation}
where $\hat{u}^{n,k}_K$ is the computed output  at the time step  $n+1$ and the iteration $k$.

One of the advantages of using ANN to build implicit linear schemes is the compatibility between this approach and the implicit step-by-step resolution. Indeed, if for $n=2$ we need to generate a randomly initial vector of  weights and the bias,    we will set $w^{n,0}=w^{n-1}$ starting from the fourth time-step. Here, $w^{n-1}$ is the vector of weights and the bias at the time-steps $n$. This is because it is expected that the behavior of the solution does not change much between two successive time-steps, and therefore  $w^n$ will be very close to $w^{n-1}$. This will reduce the number of iterations needed to obtain the final bias and weights, specially for small time-steps.

\subsection{The SSTLI method} 
The SSTLI method is an efficient version of the ESTLI method in which training of the network does not occur at each time-step, but on selected time-steps. The SST algorithm is developed to this purpose. Two criteria must be satisfied to train the network at the time step $n+1$:

\begin{enumerate}[label=\arabic*.]
\item   The performance required is not reached for the time-step $n$: It is widely expected that the approximation of $u^n$ obtained by the use of ANN is better than simply use of a semi-implicit linearization, which means that $\|\,\widetilde u^{n}-u^{n}\,\|_2 <  \, \|\, u^{n-1}-u^{n}\,\|_2$ ($\widetilde u^{n}$ is defined by \eqref{(20)} and $u^n$ is computed by the scheme \eqref{(6)}--\eqref{(EST)}). To measure the performance of the scheme, the approach proposed is to introduce a parameter $\beta>1$, and depending on the choice of this parameter, the scheme \eqref{(6)}--\eqref{(EST)} will be said to satisfy the required performance at the time-step $n$ if:
\begin{equation}
 \|\,\widetilde u^{n}-u^{n}\,\|_2 \leq \dfrac{1}{\beta} \, \|\, u^{n-1}-u^{n}\,\|_2. \label{beta}
 \end{equation}

In this case, the network is not trained for the following time-step, and we set $w^{n}=w^{n-1}$ in order to compute $\widetilde u^{n+1}$.\\

\item $\widetilde u^{n}$ is a better approximation for $u^n$ than $\overline u^{n}$: Defining $\overline u^{n}$ by, for all $K \in \T$, $\overline u_K^{n}=w_1^{n-2} u_K^{n-2}+w_2^{n-2} u_K^{n-1}+w_0^{n-2}$, the condition $\|\,\widetilde u^{n}-u^{n}\,\|_2 < \|\,\overline u^{n}-u^{n}\,\|_2$ must be satisfied to train the network at the time-step $n+1$. Otherwise, $w^{n-2}$ will be used to compute $\widetilde u^{n+1}$ and $\widetilde u^{n+2}$.
\end{enumerate} 

The detailed steps of the SST algorithm are presented in the following.  The algorithm stops when the the number of time-steps to reach the final time $T$ is exceeded.
\begin{enumerate}[label=Step \arabic*:]
\item For $n=0$ and $n=1$, compute the numerical solution using the semi-implicit scheme  \eqref{(6)}--\eqref{(7)}.
\item For $n=2$, randomly generate $w^{n,0}$, train the neural network (Fig. \ref{Fig.1}) until the error $E^{n,k}\leq \epsilon$, use  the obtained vector of weights and bias $w^n$ to compute   $\widetilde u^{n+1}$, and then solve the scheme \eqref{(6)}--\eqref{(EST)} to have $u^{n+1}$. Finally, set $n=n+1$ and go to the step 4.
\item If $\|\,\widetilde u^{n}-u^{n}\,\|_2 \geq \|\,\overline u^{n}-u^{n}\,\|_2$, set $w^{n}=w^{n-2}$ and go to the step 6. Otherwise, go to the next step.
\item If  $\|\,\widetilde u^{n}-u^{n}\,\|_2 \leq \dfrac{1}{\beta} \, \|\, u^{n-1}-u^{n}\,\|_2$, set $w^{n}=w^{n-1}$ and go to the step 9. Otherwise, set $w^{n,0}=w^{n-1}$, train the neural network (Fig. \ref{Fig.1}) until the error $E^{n,k}\leq \epsilon$ and obtain $w^{n}$, go to the next step.
\item Compute   $\widetilde u^{n+1}$ and $\overline u^{n+1}$ defined by: for all $K \in \T$, $\overline u_K^{n+1}=w_1^{n-1} u_K^{n-1}+w_2^{n-1} u_K^{n}+w_0^{n-1}$. Solve the scheme \eqref{(6)}--\eqref{(EST)} to obtain $u^{n+1}$, set $n=n+1$ and go to the step 3.
\item Compute   $\widetilde u^{n+1}$, solve the scheme \eqref{(6)}--\eqref{(EST)} to obtain $u^{n+1}$, set $n=n+1$ and go to the next step.
\item Set $w^{n}=w^{n-1}$. If $\|\,\widetilde u^{n}-u^{n}\,\|_2 \leq \dfrac{1}{\beta} \, \|\, u^{n-1}-u^{n}\,\|_2$  go to the step 9. Otherwise go to the next step.
\item  Compute   $\widetilde u^{n+1}$, solve the scheme \eqref{(6)}--\eqref{(EST)} to obtain $u^{n+1}$, and set $\overline u^{n+1}=\widetilde u^{n+1} $.  Then, set $w^{n,0}=w^{n-1}$, train the neural network (Fig. \ref{Fig.1}) until the error $E^{n,k}\leq \epsilon$ and obtain a new $w^{n}$. Compute a new $\widetilde u^{n+1}$, set $n=n+1$  and go to the step 3.
\item Compute   $\widetilde u^{n+1}$, solve the scheme \eqref{(6)}--\eqref{(EST)} to obtain $u^{n+1}$. Set $n=n+1$ and go to the step 4.
 
\end{enumerate}

 \section{Numerical experiments}
 
  This section mainly deals with  the comparison between the Euler semi-implicit method, the ESTLI method, and the SSTLI method. In all numerical tests presented in this section, target outputs are normalized to the range of $[-1,1]$, and the mean squared error  is required to be less than $\epsilon=10^{-3}$.  The Levenberg-Marquard  method is chosen as learning algorithm to train the single-layer neural network, and a minimum performance gradient of $10^{-10}$ is required for more precision. Finally, for the parameter $\beta$ defined in \eqref{beta}, we take $\beta=10$.
 
 \subsection{Keller-Segel model for embryonic pattern formation}

In this test, we consider the system \eqref{(2)}--\eqref{(3)}. The spatial  domain is $\Omega=(-7/2,7/2)\times(-35,35) $ with a uniform mesh of $12250$ control volumes, and the final time is $T=150$. We adopt the following initial condition:
\begin{equation}
  u_0(x) = \left\{\begin{array}{ll}
1+\alpha(x) &\quad\mbox{if } x\in (-7/2,7/2)\times(-1,1),\\ 
 1 &\quad\mbox{otherwise, with  } x \in \Omega,
  \end{array}\right. \\
  \label{(init)}
  \end{equation}
  where $\alpha(x)$ is a positive perturbation  which is constant on each control volume. It is equal on each cell of the mesh  to the average of ten uniformly distributed random values in $[0,1]$. For the model parameters, we take: $D_u=0.25$ and $\chi=2$.
  
\begin{figure}[h!]
\subfigure{\includegraphics[width=6.cm]{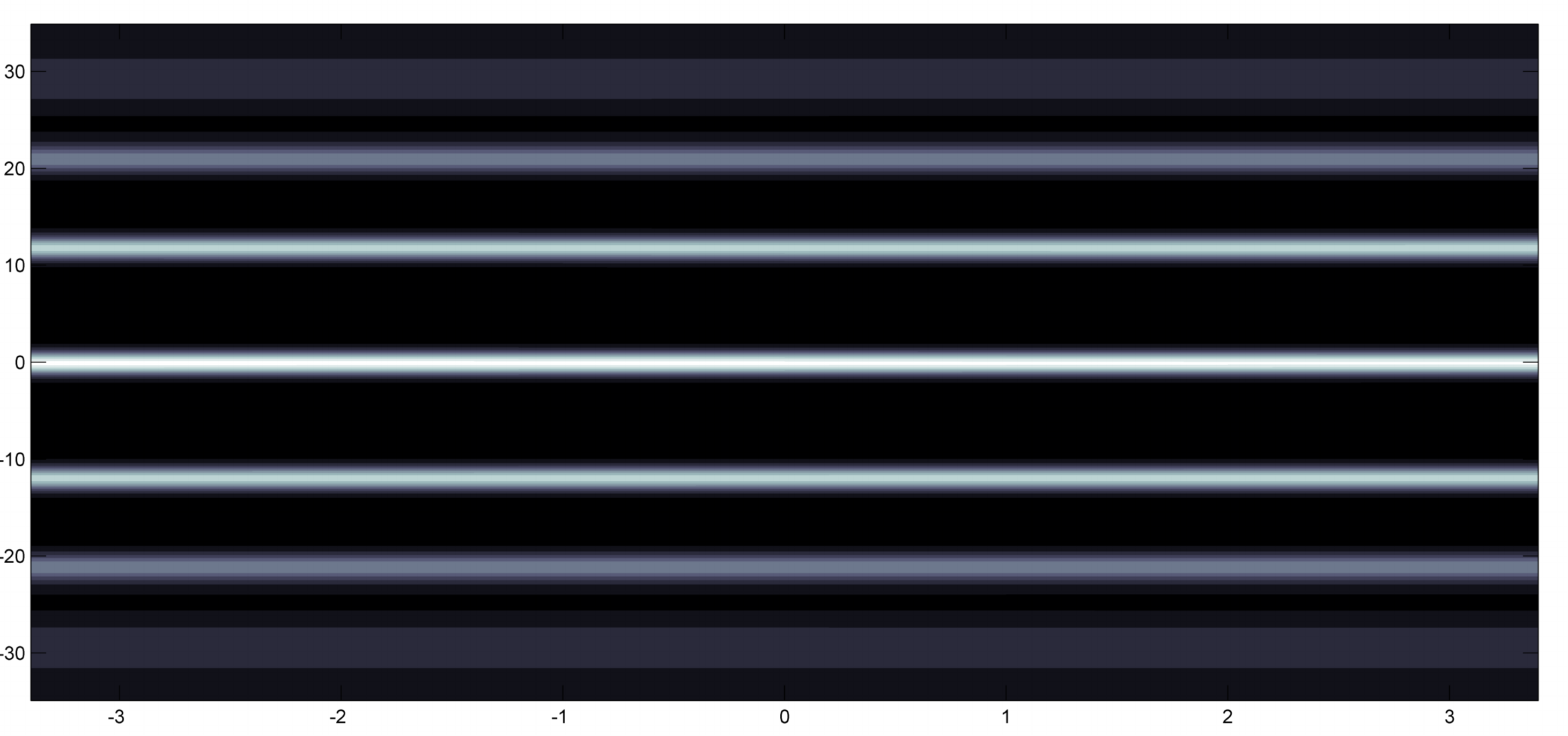}}
\subfigure{\includegraphics[width=6.cm]{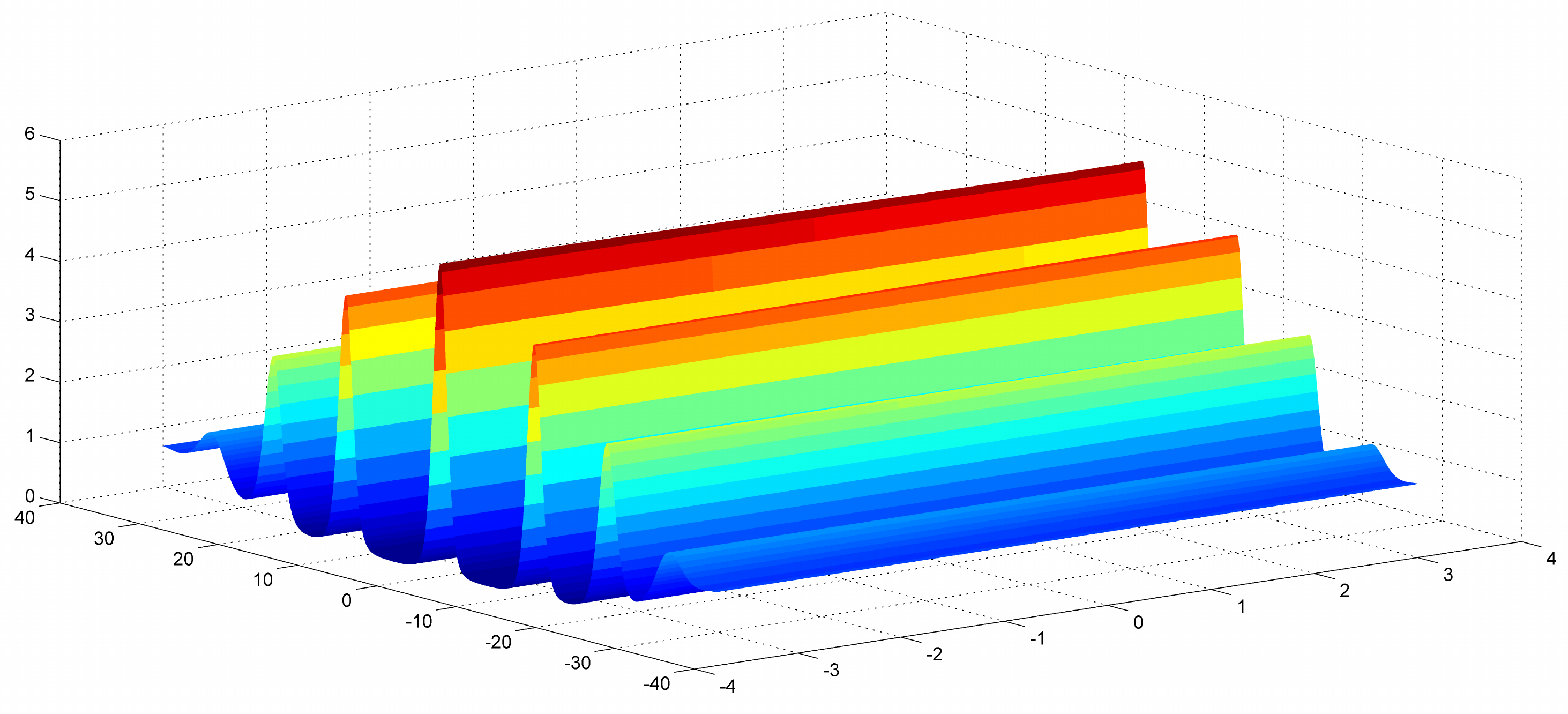}}
\caption{Solution $(u)$ of  \eqref{(2)}--\eqref{(3)} at final time  computed via the SSTLI method with $\Delta t=10^{-2}$ }
\label{Fig.3}
\end{figure}

\begin{figure}[h!]
\subfigure{\includegraphics[width=6.cm]{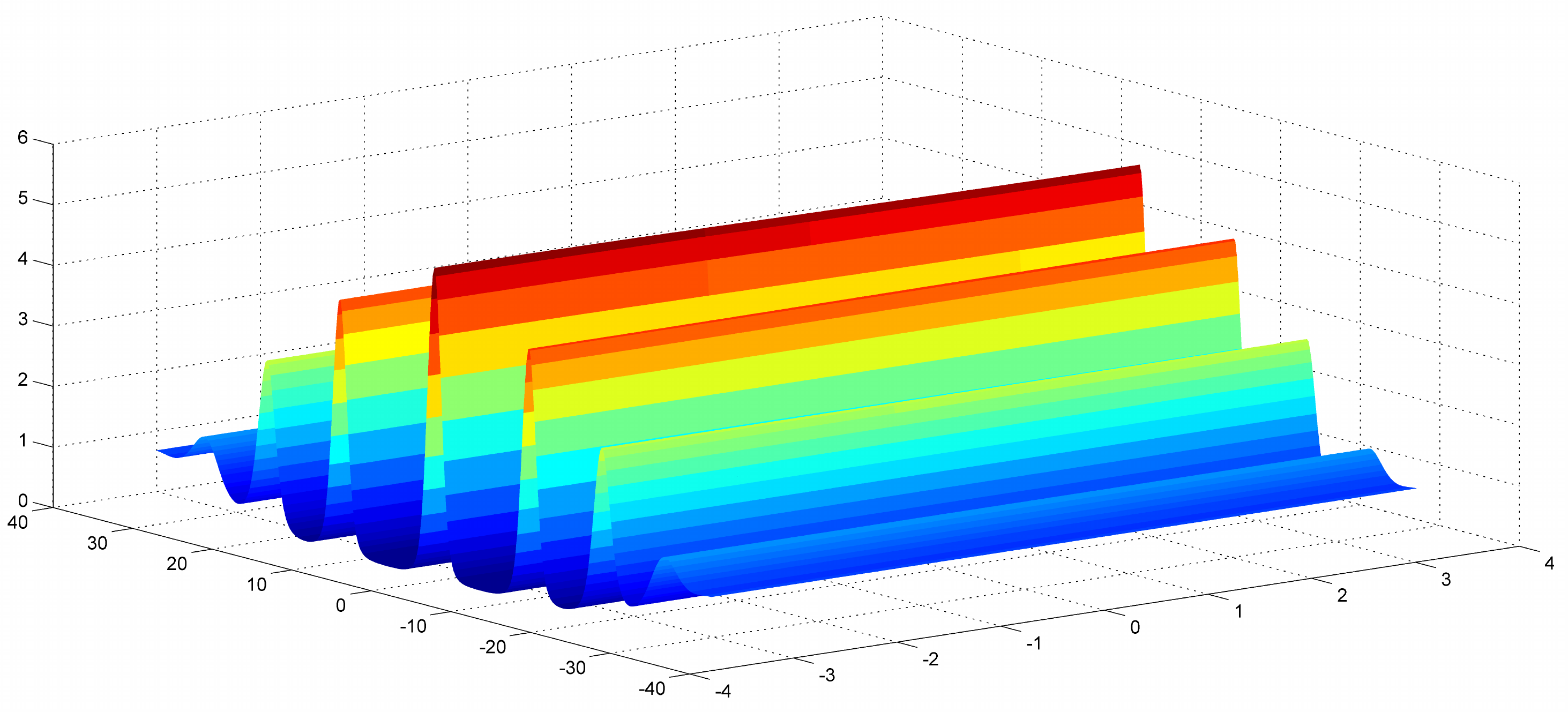}}
\subfigure{\includegraphics[width=6.cm]{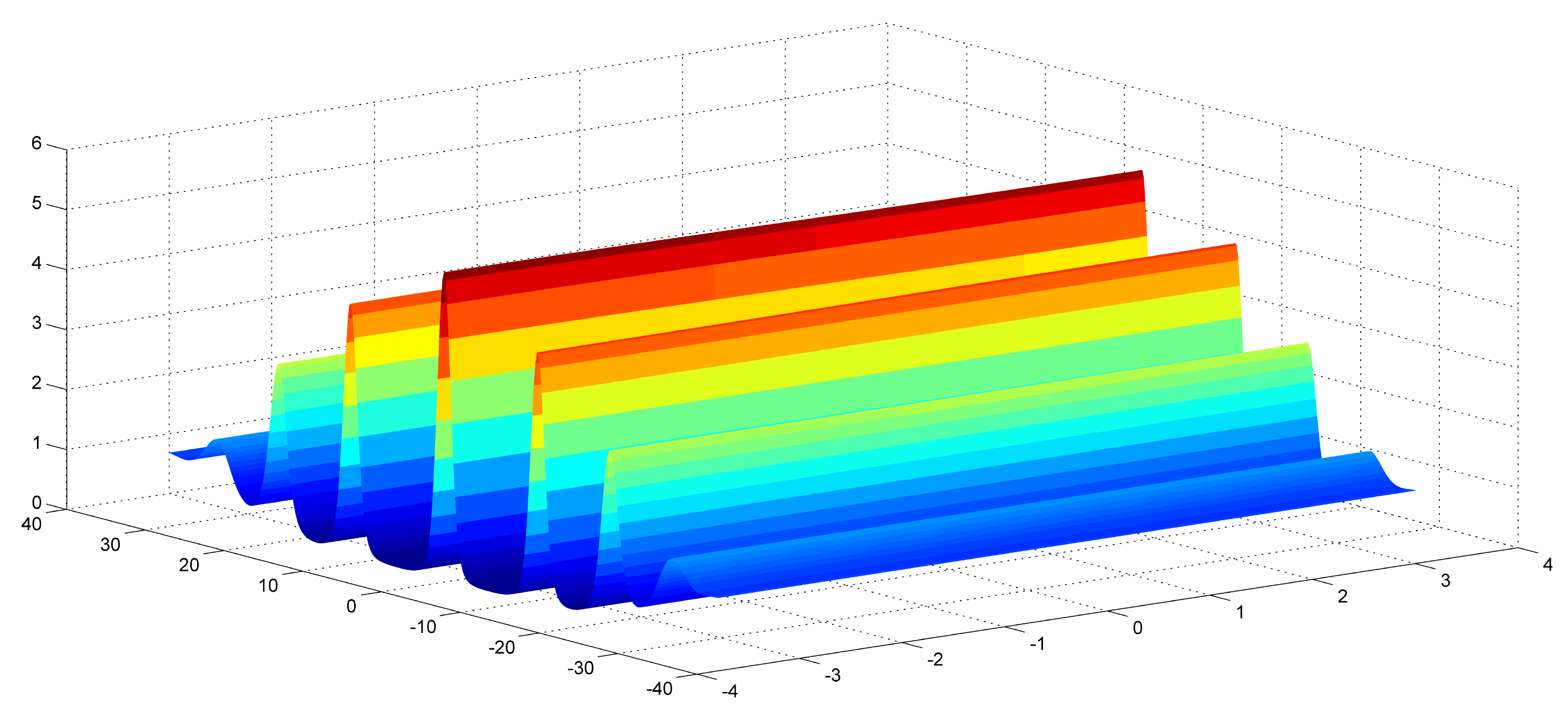}}
\caption{Solution $(u)$ of  \eqref{(2)}--\eqref{(3)} at final time  computed via the SSTLI method (left) and the ESTLI method (right) with $\Delta t=5$ }
\label{Fig.4}
\end{figure}

  The numerical cell density $u$ obtained by the SSTLI method with $\Delta t=10^{-2}$ is presented in Fig. \ref{Fig.3}.  We can observe from this  figure  the ability of the studied model  to generate stripe patterns. We can see also from the three-dimensional plot on the right that the obtained numerical solution is positive and is oscillation free. In  Fig. \ref{Fig.4}, we present three-dimensional plots of the cell density using SSTLI and ESTLI methods with a very large time-step size $\Delta t=5$. As we can see, both numerical solutions are free from negative values or any instabilities.

\begin{figure}[h!]
\begin{center}
\includegraphics[width=0.75\textwidth]{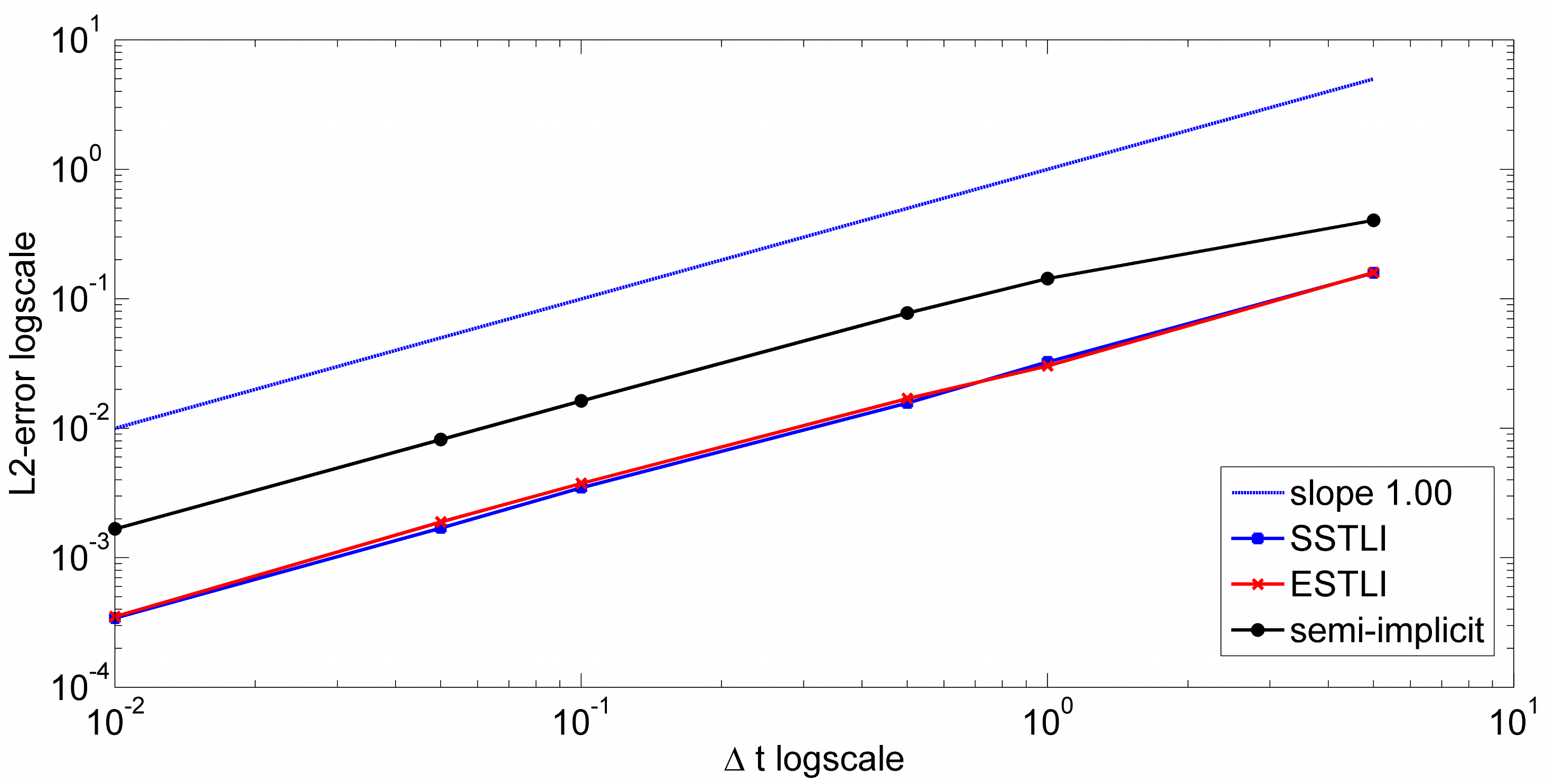}
\caption{Convergence speed of the SSTLI, ESTLI, and Euler semi-implicit methods applied to the the model \eqref{(2)}--\eqref{(3)} }
\label{Fig.5}
\end{center}
\end{figure}

 In order to compare the accuracy of the three approaches presented in this paper, and since the exact  solution of the Keller-Segel system \eqref{(2)}--\eqref{(3)} is unavailable,  we use a reference solution  computed by the  corrected decoupled scheme  \cite{Akhmouch2} on a very fine  time-stepping $\Delta t=10^{-3}$. This reference solution is  used to compute the relative $L^2$-errors at final time (see Table \ref{tab.1}) for several  time-step sizes.  Computational cost (CPU) of the  semi-implicit scheme \eqref{(6)}--\eqref{(8)}, and the increase in computational cost (denoted as $\Gamma$) of SSTLI and ESTLI methods over the semi-implicit scheme are also presented in Table \ref{tab.1}. From this table, we can observe that  SSTLI method is slightly more accurate than ESTLI method except for $\Delta t=1$. Moreover, both methods are in general four to five times more accurate than  Euler semi-implicit scheme.  Concerning the computational time, the presented results show that ESTLI method is considerably more time-consuming than  semi-implicit
scheme. On the contrary, the increase in computational cost over  semi-implicit scheme is generally weak in the case of SSTLI method, specially for the three finest time-steps (less than $15 \%$). Finally, in  Fig. \ref{Fig.5}, it can be seen that the three numerical approaches are first-order accurate.
\begin{table}[h!]
\caption{Computational time and relative $L^2$-errors  obtained for  the solution $(u)$ of  \eqref{(2)}--\eqref{(3)} using SSTLI, ESTLI, and Euler semi-implicit methods}
\label{tab.1}       
\begin{tabular}{lllllll}
\hline\noalign{\smallskip}
{$\Delta t$ } & $L^2$-error  & {$\Gamma $} ($\%$) & $L^2$-error ESTLI & {$\Gamma $} ($\%$) & $L^2$-error  & CPU (s)  \\
 &  SSTLI &  &  ESTLI &  &  semi-implicit &   \\
\noalign{\smallskip}\hline\noalign{\smallskip}
$5$ & $1.585 \times 10^{-1}$ & $168.87$  & $ 1.596\times 10^{-1}$ &  $ 171.01$ &  $4.042 \times 10^{-1}$ &  $7.2$  \\
$1$ & $ 3.257\times 10^{-2}$ & $ 28.36$  & $ 3.036 \times 10^{-2}$ &   $184.76$ & $1.435 \times 10^{-1}$ &  $32.4$  \\
$5.10^{-1}$ & $ 1.569 \times 10^{-2}$ & $ 19.59$  & $ 1.698 \times 10^{-2}$ &  $ 182.81$ & $7.767 \times 10^{-2}$ &  $65.0$  \\
$10^{-1}$ & $3.478 \times 10^{-3}$ & $ 13.97$  & $3.760 \times 10^{-3}$ &  $ 154.45$ & $1.630 \times 10^{-2}$ &  $315.8$  \\
$5.10^{-2}$ & $1.701 \times 10^{-3}$ & $ 13.17$  & $ 1.893 \times 10^{-3}$ &  $ 118.63$ & $8.205 \times 10^{-3} $ &  $637.7$  \\
$10^{-2}$ & $3.424 \times 10^{-4}$ & $ 12.58$  & $ 3.512 \times 10^{-4}$ &  $ 84.29$ & $1.672 \times 10^{-3}$ &  $3186.1$  \\
\noalign{\smallskip}\hline
\end{tabular}
\end{table}

\subsection{Keller-Segel model with quadratic growth}

In this subsection, we are concerned with the system  \eqref{(4)} endowed with zero-flux boundary conditions. Herein, the computational domain is $\Omega=(-8,8)^2$  with a uniform mesh grid  $100 \times 100$., $T=30$, $D_u=0.0625$, $\chi=6$, $\lambda=16$ and $f(u)=2u(1-u)$. The following initial conditions are considered:
\begin{equation*}
  u(x,0) = \left\{\begin{array}{ll}
1+\alpha(x) &\quad\mbox{if } \|x\|_2<0.7,\\ 
 1 &\quad\mbox{otherwise},
  \end{array}\right. \\
\end{equation*}
and  $c(x,0)=1/32$. The  perturbation $\alpha(x)$ is defined as in \eqref{(init)}. The scheme adopted is:
\begin{align}
& \m(K)\frac{u^{n+1}_K-u^n_K}{\Delta t}
  - D_u \sum_{\sigma\in\E_K}\tau_\sigma Du_{K,\sigma}^{n+1} \notag
  \\
  &+\chi \sum_{\substack{\sigma\in\E_{K}\\ \sigma=K|L}}\tau_\sigma\left(S\left( Dc_{K,\sigma}^{n+1}\right)u_{K}^{n+1}-S\left( -Dc_{K,\sigma}^{n+1}\right)u_{L}^{n+1}\right)-2\m(K)\,\widetilde u_K^{n+1}\left(1-u_K^{n+1}\right)=0, \label{x1}
 \\
  &\m(K)\frac{c^{n+1}_K-c^n_K}{\Delta t} -\sum_{\sigma\in\E_K}\tau_\sigma\, Dc^{n+1}_{K,\sigma}
  +\lambda \m(K)\, c_K^{n+1}= \m(K)\, \widetilde u_K^{n+1}, \label{x2}
\end{align}
when using the semi-implicit method, $\widetilde u_K^{n+1}$ is replaced by $u_K^{n}$, for ESTLI method $\widetilde u_K^{n+1}$ is computed by an ANN similar to that of Figs. \ref{Fig.1} and \ref{Fig.2}. In the case of SSTLI method, SST algorithm is used.

\begin{figure}[h!]
\subfigure{\includegraphics[width=6.cm]{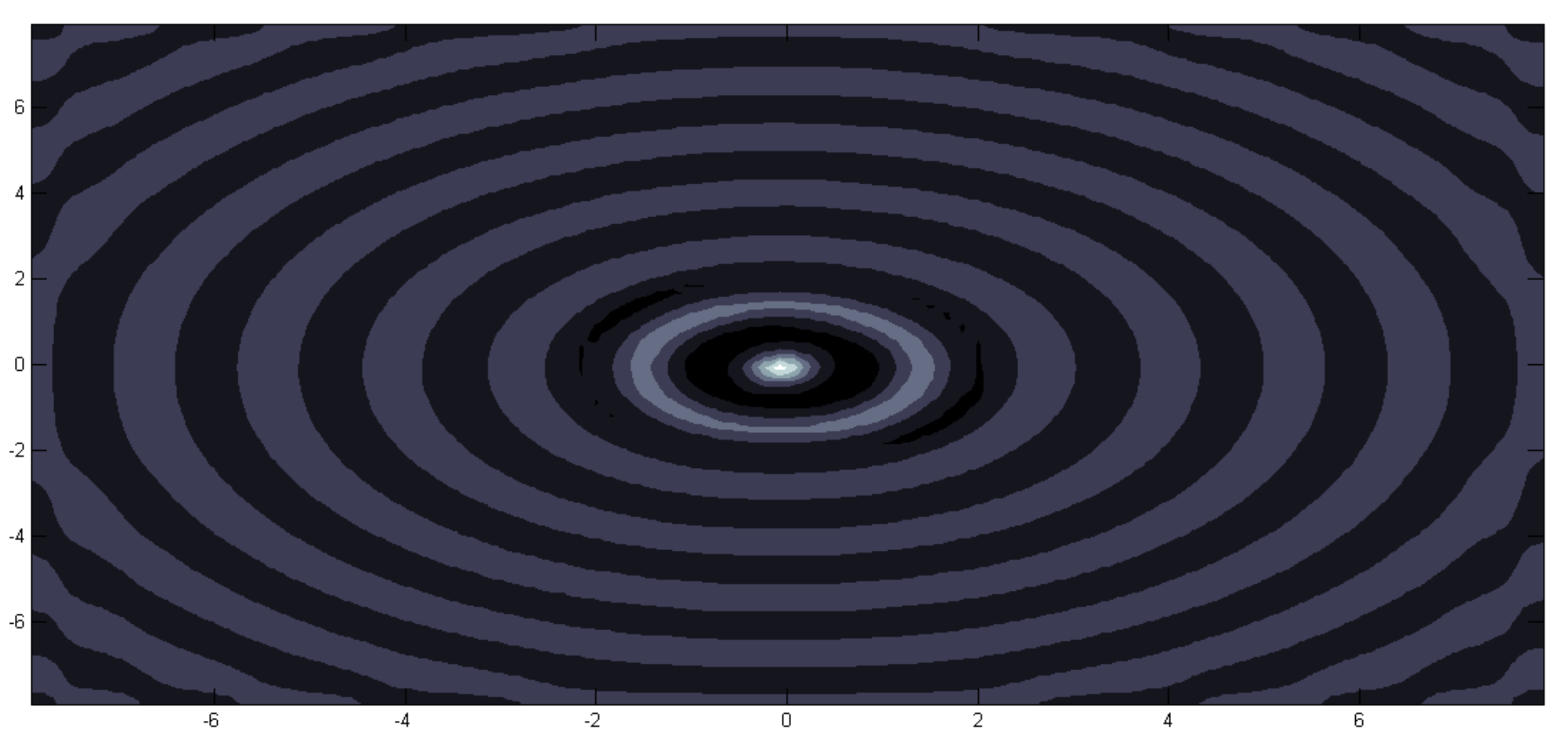}}
\subfigure{\includegraphics[width=6.cm]{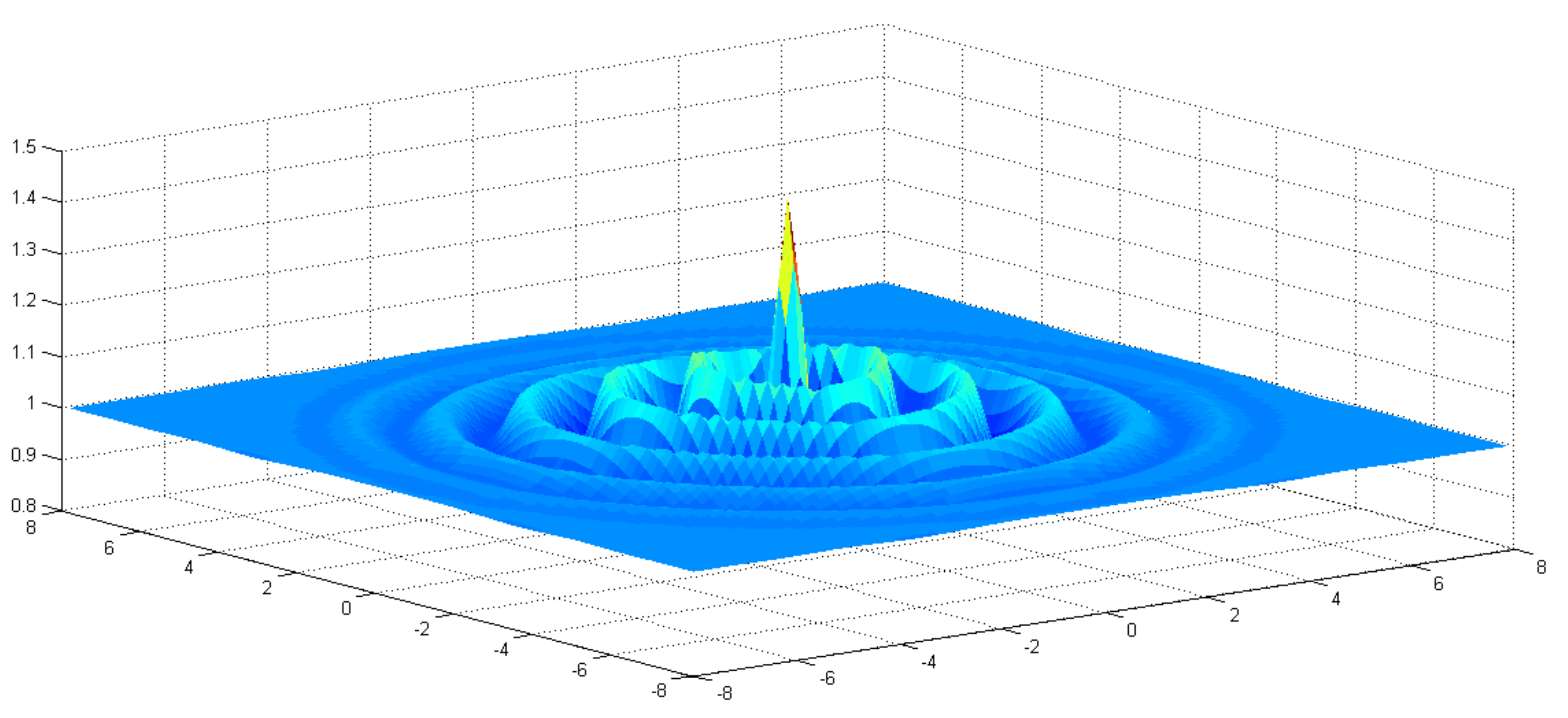}}
\caption{Solution $(u)$ of  \eqref{(4)} with quadratic growth at final time  computed via the SSTLI method with $\Delta t=5.10^{-3}$ }
\label{Fig.6}
\end{figure}

\begin{figure}[h!]
\subfigure{\includegraphics[width=6.cm]{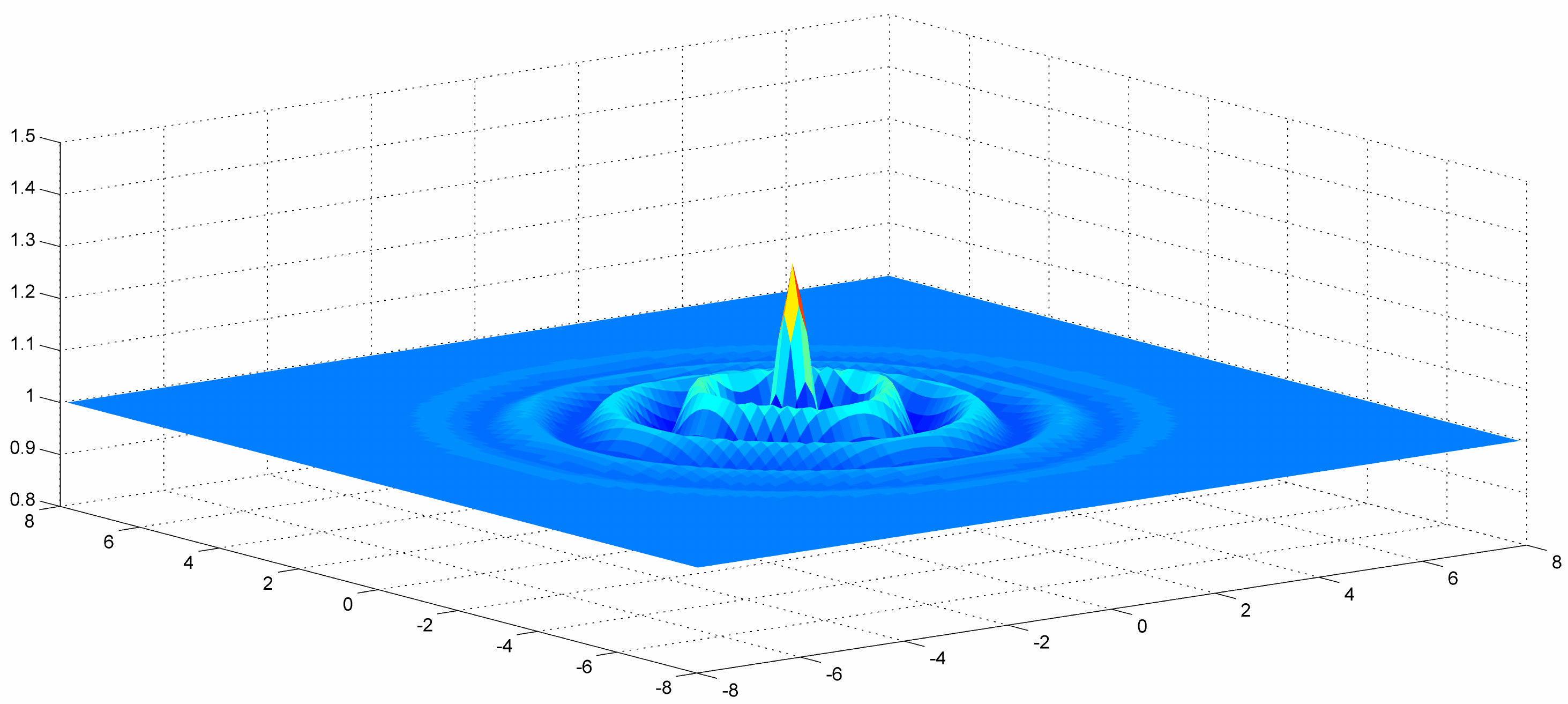}}
\subfigure{\includegraphics[width=6.cm]{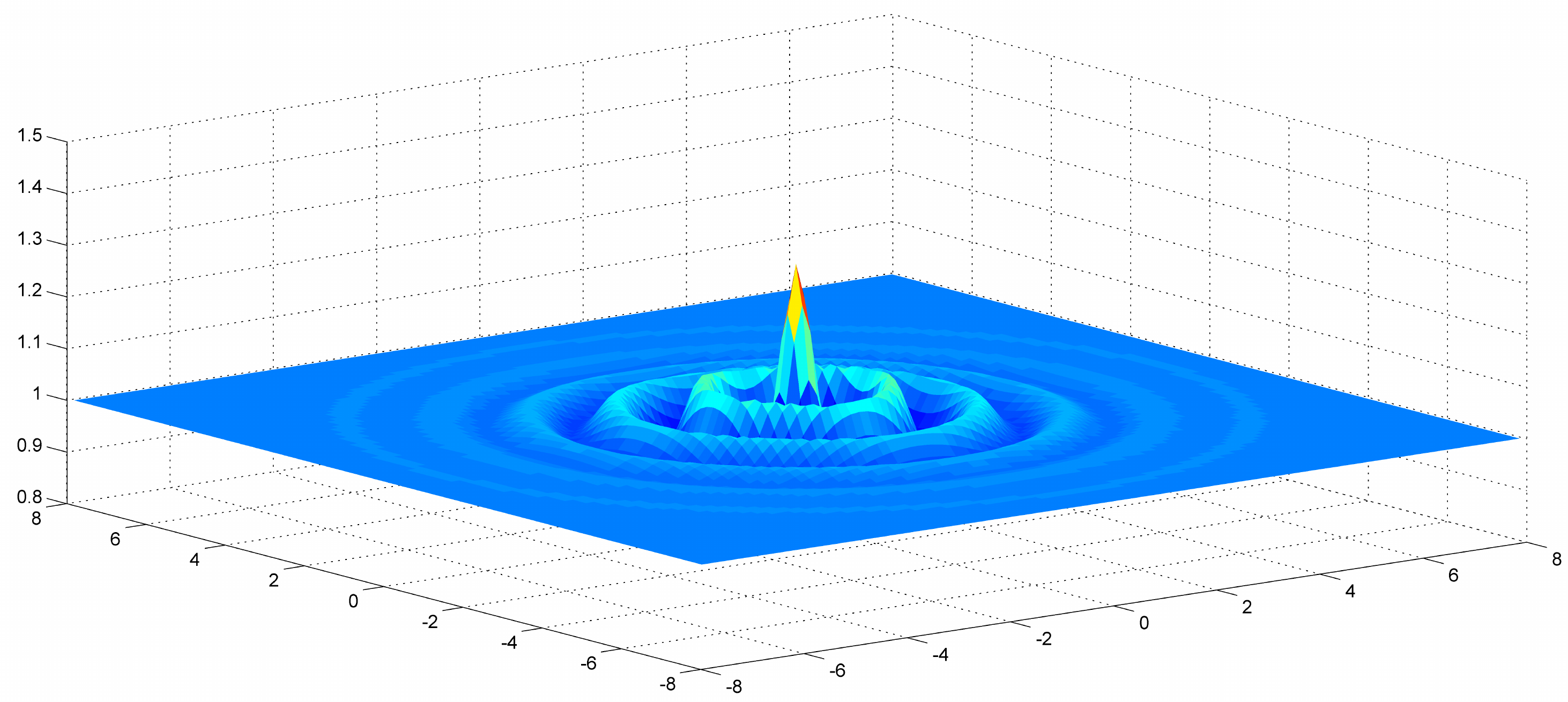}}
\caption{Solution $(u)$ of  \eqref{(4)} with quadratic growth at final time  computed via the SSTLI method (left) and the ESTLI method (right) with $\Delta t=1$ }
\label{Fig.7}
\end{figure}

The numerical cell density $u$ computed via SSTLI method with $\Delta t=5.10^{-3}$  is plotted in Fig. \ref{Fig.6}. The continuous rings observed in the figure are similar to the patterns formed by \textit{Salmonella typhimurium} bacteria \cite{Woodward1}. Numerical solutions computed by the SSTLI and ESTLI methods with $\Delta t=1$ are shown in Fig. \ref{Fig.7}.  We can see that even for a large time-step, the stability and the positivity of the solution are assured. 

\begin{figure}[h!]
\begin{center}
\includegraphics[width=0.75\textwidth]{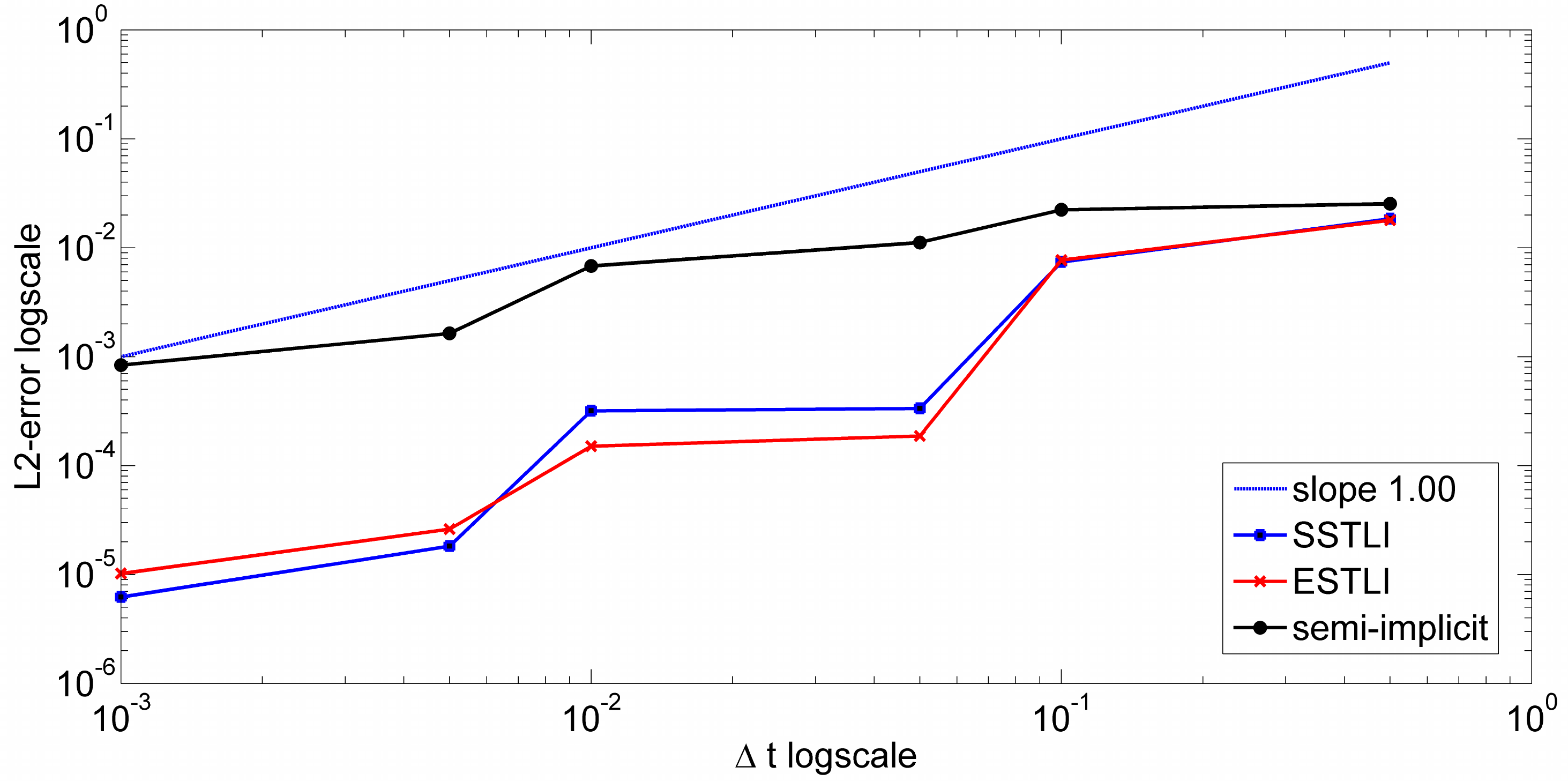}
\caption{Convergence speed of the SSTLI, ESTLI, and Euler semi-implicit methods applied to the the model \eqref{(4)} with quadratic growth }
\label{Fig.8}
\end{center}
\end{figure}

The reference solution used is computed using the corrected decoupled scheme with $\Delta t=10^{-4}$. For a sequence of decreasing time stepping, relative $L^2$-errors and  computational costs of the studied numerical methods are presented in Table \ref{tab:2}, convergence speeds are presented in Fig. \ref{Fig.8}. Table \ref{tab:2} shows that SSTLI and ESTLI methods are much more accurate than semi-implicit method. The obtained results do not allow us to conclude which of the two proposed methods is more accurate, but it can be observed from the table that SSTLI give better results for the two smallest time-steps. In terms of computational cost, the increase in computational time  over  semi-implicit scheme is very marginal for the three most time-consuming tests when using SSTLI method (less than $10 \%$), which is not the case for the ESTLI method.
\begin{table}[h!]
\caption{Computational time and relative $L^2$-errors  obtained for  the solution $(u)$ of  \eqref{(4)} with quadratic growth using SSTLI, ESTLI, and Euler semi-implicit methods}
\label{tab:2}       
\begin{tabular}{lllllll}
\hline\noalign{\smallskip}
{$\Delta t$ } & $L^2$-error  & {$\Gamma $} ($\%$) & $L^2$-error  & {$\Gamma $} ($\%$) & $L^2$-error  & CPU (s)  \\
 &  SSTLI &  &  ESTLI &  &  semi-implicit &   \\
\noalign{\smallskip}\hline\noalign{\smallskip}
$1$ & $1.850 \times 10^{-2}$ & $228.52$  & $ 1.790\times 10^{-2}$ &  $ 273.21$ &  $2.536\times 10^{-2}$ &  $5.6$  \\
$5.10^{-1}$ & $7.425 \times 10^{-3}$ & $65.64$  & $ 7.703\times 10^{-3}$ &  $ 229.49$ &  $2.234\times 10^{-2}$ &  $11.6$  \\
$10^{-1}$ & $ 3.346\times 10^{-4}$ & $ 20.02$  & $ 1.873 \times 10^{-4}$ &   $165.30$ & $1.119\times 10^{-2}$ &  $55.3$  \\
$5.10^{-2}$ & $ 3.185 \times 10^{-4}$ & $ 8.78$  & $ 1.508 \times 10^{-4}$ &  $ 140.80$ & $6.806\times 10^{-3}$ &  $110.1$  \\
$10^{-2}$ & $1.823 \times 10^{-5}$ & $ 7.72$  & $2.614 \times 10^{-5}$ &  $ 65.17$ & $1.636\times 10^{-3}$ &  $538.9$  \\
$5.10^{-3}$ & $6.226 \times 10^{-6}$ & $ 6.74$  & $ 1.022 \times 10^{-5}$ &  $ 51.21$ & $8.385\times 10^{-4} $ &  $1090.8$  \\
\noalign{\smallskip}\hline
\end{tabular}
\end{table}

\subsection{Keller-Segel model with volume-filling }
Here, we focus on the model \eqref{(5)} with zero-flux boundary conditions. The computational domain and the initial chemoattractant concentration are taken as in the previous subsection, and we prescribe the following initial cell density
\begin{equation*}
  u(x,0) = \left\{\begin{array}{ll}
\alpha(x) &\quad\mbox{if } \|x\|_2<0.7,\\ 
 0 &\quad\mbox{otherwise},
  \end{array}\right. \\
\end{equation*}
$\alpha(x)$ is defined as in \eqref{(init)}. Besides, we set $D_u=0.1$ and $\chi(u)=10(1-u)$. The following scheme is used:
\begin{align*}
& \m(K)\frac{u^{n+1}_K-u^n_K}{\Delta t}
  - D_u \sum_{\sigma\in\E_K}\tau_\sigma Du_{K,\sigma}^{n+1} \notag
  \\
  &+ \sum_{\substack{\sigma\in\E_{K}\\ \sigma=K|L}}\tau_\sigma\left(S\left(\scalebox{1.2}{$\widetilde \chi$}^{n+1}_{K,\sigma} Dc_{K,\sigma}^{n+1}\right)u_{K}^{n+1}-S\left( -\scalebox{1.2}{$\widetilde \chi$}^{n+1}_{K,\sigma} Dc_{K,\sigma}^{n+1}\right)u_{L}^{n+1}\right)=0, 
 \\
  &\m(K)\frac{c^{n+1}_K-c^n_K}{\Delta t} -\sum_{\sigma\in\E_K}\tau_\sigma\, Dc^{n+1}_{K,\sigma}
  + \m(K)\, c_K^{n+1}= \m(K)\, \widetilde u_K^{n+1}, 
\end{align*}
 where the value of $\widetilde u_K^{n+1}$ is determined as in \eqref{x1}--\eqref{x2}, and $\scalebox{1.2}{$\widetilde \chi$}^{n+1}_{K,\sigma}=\dfrac{\chi \left(\widetilde u_K^{n+1} \right)+\chi \left(\widetilde u_L^{n+1} \right)}{2}$ when $\sigma=K|L$.
 
\begin{figure}[h!]
\subfigure{\includegraphics[width=6.cm]{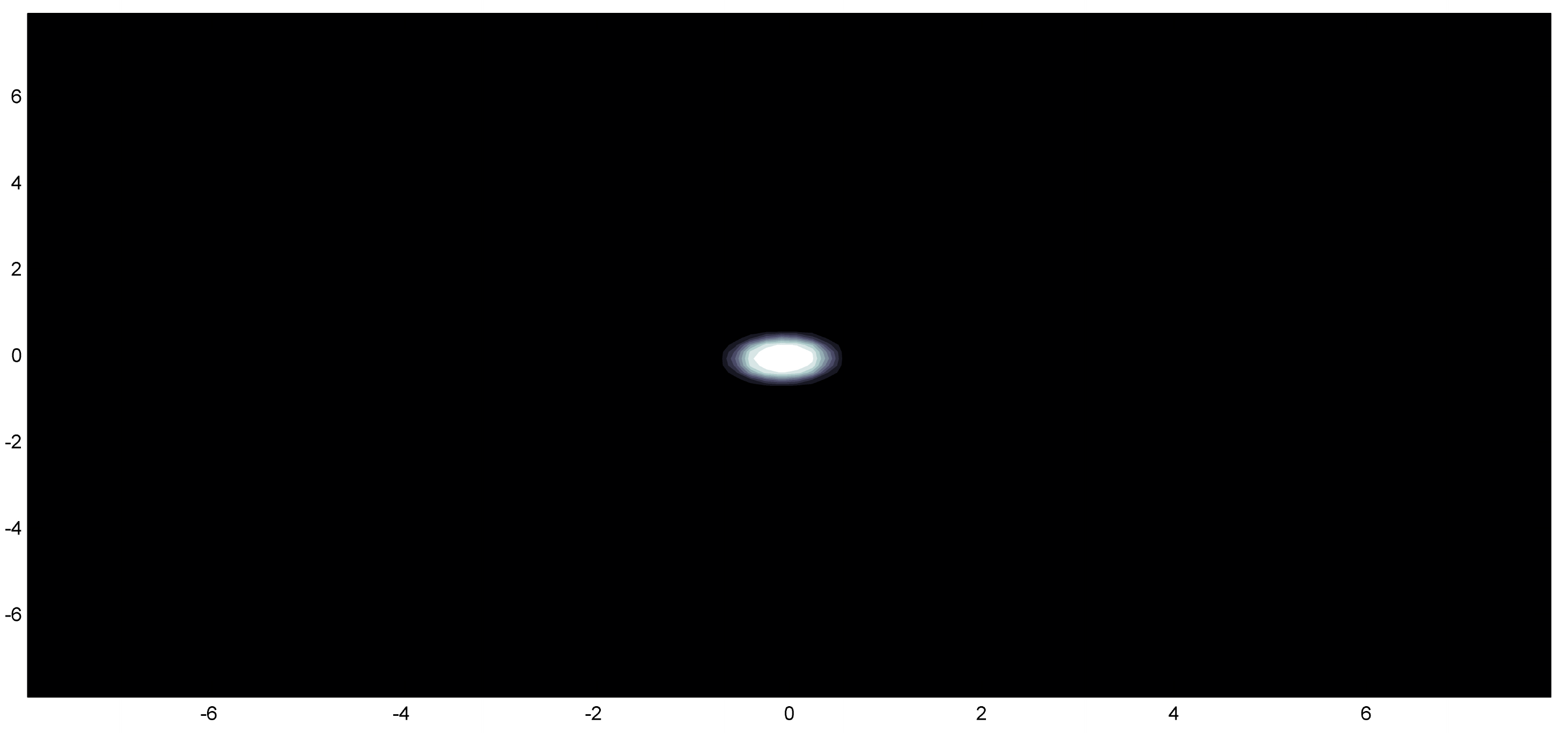}}
\subfigure{\includegraphics[width=6.cm]{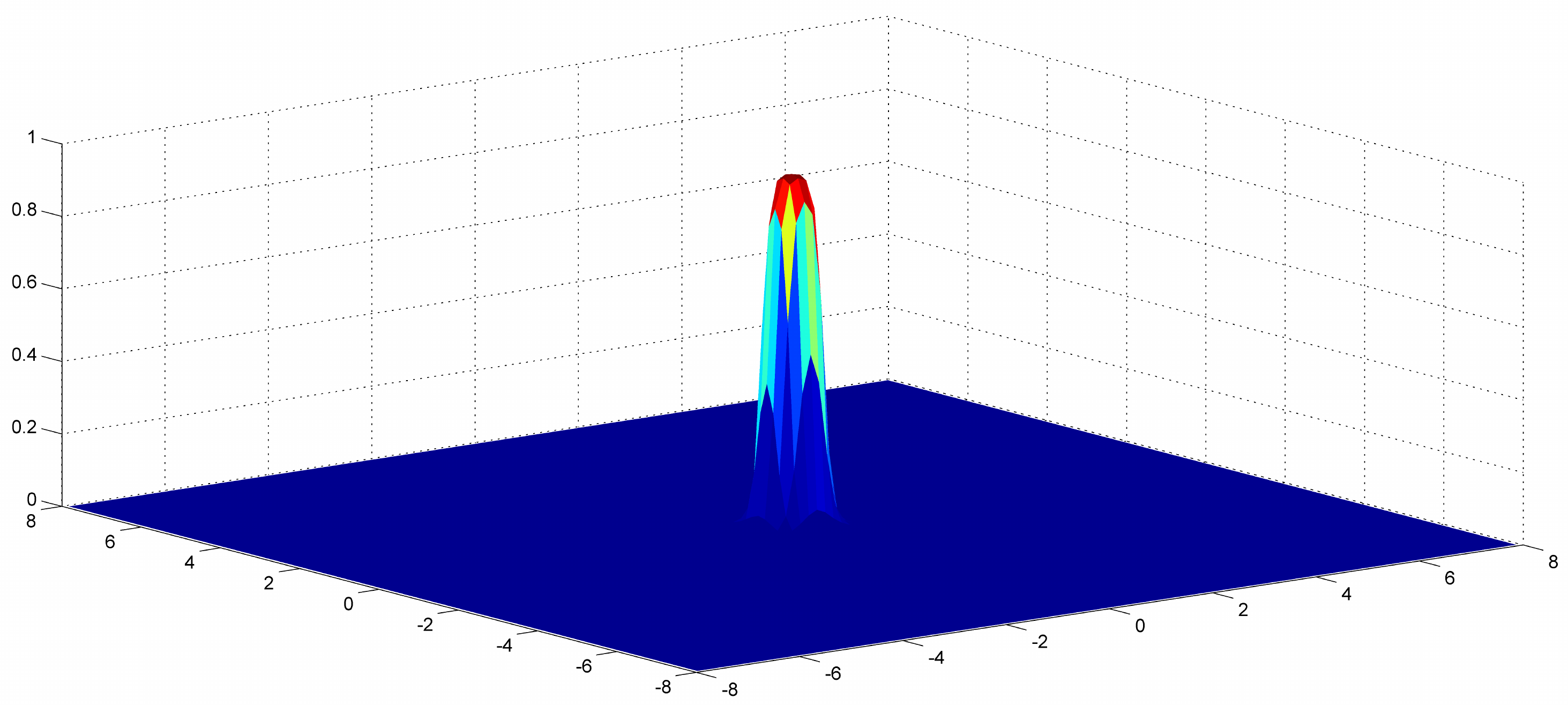}}
\caption{Solution $(u)$ of  \eqref{(5)} at final time  computed via the SSTLI method with $\Delta t=10^{-4}$ }
\label{Fig.9}
\end{figure}

 Numerical solution computed by SSTLI method with $\Delta t=10^{-4}$  is plotted in Fig. \ref{Fig.9}. We observe from the figure that the maximal  density of cells does not exceed $1$ even if the solution concentrates at the center, which is consistent with our choice of the function $\chi$. 

\begin{figure}[h!]
\begin{center}
\includegraphics[width=0.75\textwidth]{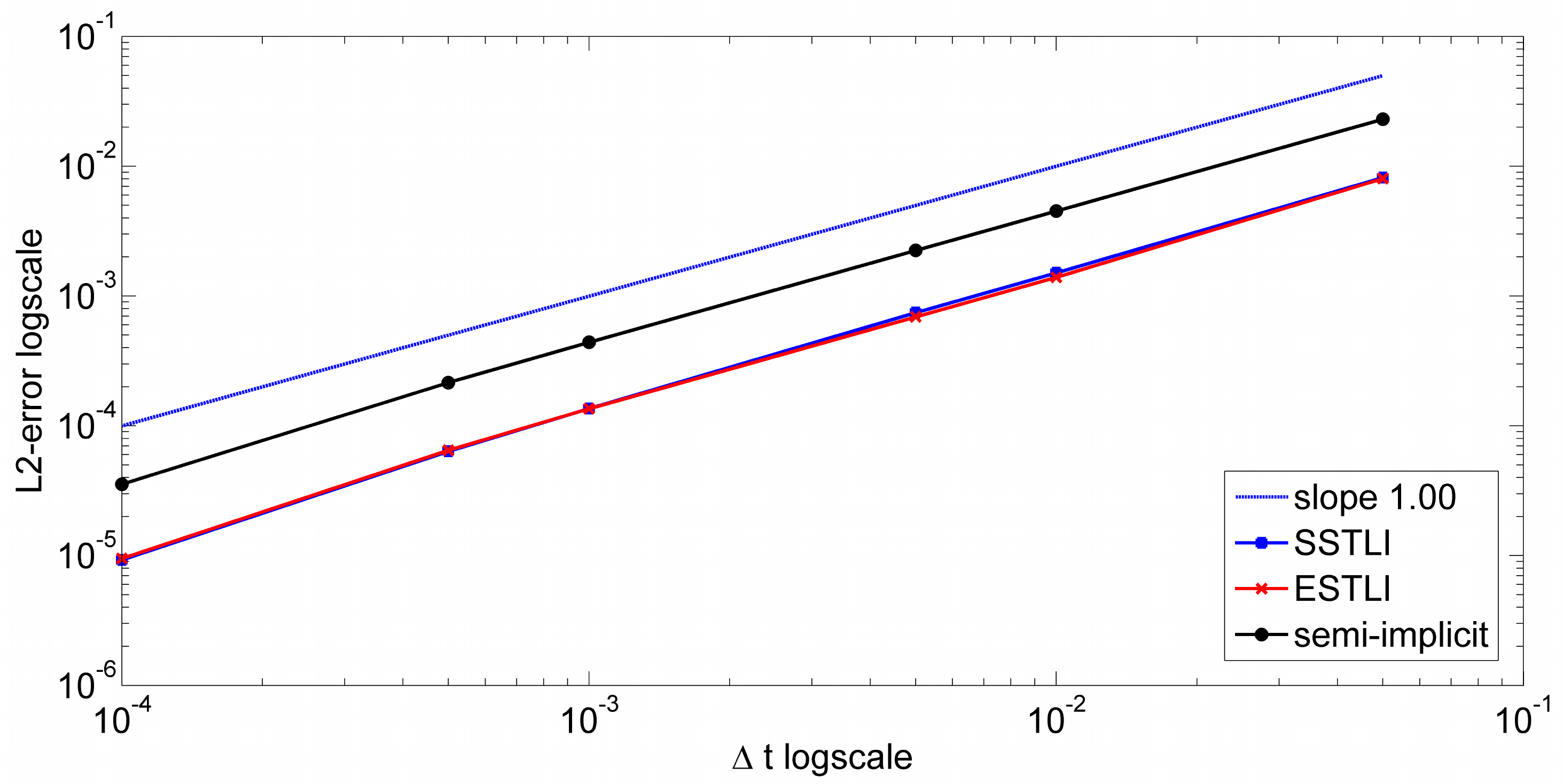}
\caption{Convergence speed of the SSTLI, ESTLI, and Euler semi-implicit methods applied to the the model \eqref{(5)} }
\label{Fig.10}
\end{center}
\end{figure}

To compare relative errors for the three numerical methods studied in this paper, the reference solution of this test is computed using the SSTLI method with $\Delta t=10^{-5}$. Table \ref{tab:3} shows that  relative $L^2$-errors for SSTLI and ESTLI methods are very close with a slight advantage for SSTLI method in the case of small time-steps. The two methods are about three to four times more accurate than semi-implicit scheme, and the   three numerical approaches are first-order accurate (see  Fig. \ref{Fig.10}). Table \ref{tab:3} also shows that, except for $\Delta t=5.10^{-2}$, the increase in computational cost over  semi-implicit scheme is  weak to insignificant for SSTLI method. In terms of $L^\infty$-error, we observe in Table \ref{tab:4} that  SSTLI and ESTLI methods are four to six times more accurate than semi-implicit scheme.

\begin{table}[h!]
\caption{Computational time and relative $L^2$-errors  obtained for  the solution $(u)$ of  \eqref{(5)} using SSTLI, ESTLI, and Euler semi-implicit methods}
\label{tab:3}       
\begin{tabular}{lllllll}
\hline\noalign{\smallskip}
{$\Delta t$ } & $L^2$-error  & {$\Gamma $} ($\%$) & $L^2$-error  & {$\Gamma $} ($\%$) & $L^2$-error  & CPU (s)  \\
 &  SSTLI &  &  ESTLI &  &  semi-implicit &   \\
\noalign{\smallskip}\hline\noalign{\smallskip}
$5.10^{-2}$ & $8.174 \times 10^{-3}$ & $70.92$  & $ 8.038\times 10^{-3}$ &  $ 149.66$ &  $2.309\times 10^{-2}$ &  $5.0$  \\
$10^{-2}$ & $ 1.510\times 10^{-3}$ & $ 13.68$  & $ 1.395 \times 10^{-3}$ &   $164.82$ & $4.518\times 10^{-3}$ &  $20.1$  \\
$5.10^{-3}$ & $ 7.458 \times 10^{-4}$ & $ 6.77$  & $ 6.917 \times 10^{-4}$ &  $ 129.15$ & $2.248\times 10^{-3}$ &  $39.4$  \\
$10^{-3}$ & $1.359 \times 10^{-4}$ & $ 4.14$  & $1.356 \times 10^{-4}$ &  $ 74.07$ & $4.409\times 10^{-4}$ &  $193.2$  \\
$5.10^{-4}$ & $6.370 \times 10^{-5}$ & $ 3.43$  & $ 6.487 \times 10^{-5}$ &  $ 54.27$ & $2.156\times 10^{-4} $ &  $390.5$  \\
$10^{-4}$ & $9.265 \times 10^{-6}$ & $ 2.37$  & $ 9.504 \times 10^{-6}$ &  $ 38.74$ & $3.553\times 10^{-5}$ &  $1966.8$  \\
\noalign{\smallskip}\hline
\end{tabular}
\end{table}

\begin{table}[h!]
\centering
\caption{Relative $L^\infty$-errors  obtained for  the solution $(u)$ of  \eqref{(5)} using SSTLI, ESTLI, and Euler semi-implicit methods}
\label{tab:4}       
\begin{tabular}{llll}
\hline\noalign{\smallskip}
{$\Delta t$ } & $L^\infty$-error   & $L^\infty$-error  & $L^\infty$-error    \\
 &  SSTLI &  ESTLI &  semi-implicit    \\
\noalign{\smallskip}\hline\noalign{\smallskip}
$5.10^{-2}$ & $7.754 \times 10^{-3}$ &  $ 7.682\times 10^{-3}$ &    $3.314\times 10^{-2}$   \\
$10^{-2}$ & $ 1.479\times 10^{-3}$ &  $ 1.370 \times 10^{-2}$ &    $6.770\times 10^{-3}$   \\
$5.10^{-3}$ & $ 7.411 \times 10^{-4}$ &  $ 6.837 \times 10^{-4}$ &   $3.383\times 10^{-3}$   \\
$10^{-3}$ & $1.412 \times 10^{-4}$ & $1.345 \times 10^{-4}$ &   $6.664\times 10^{-4}$  \\
$5.10^{-4}$ & $6.370 \times 10^{-5}$   & $ 6.415 \times 10^{-5}$ &   $3.267\times 10^{-4} $  \\
$10^{-4}$ & $9.689 \times 10^{-6}$   & $ 9.973 \times 10^{-6}$ &   $5.480\times 10^{-5}$   \\
\noalign{\smallskip}\hline
\end{tabular}
\end{table}

\subsection{Keller-Segel model with cubic growth}
The purpose of this test is to investigate the ability of SSTLI and ESTLI methods to capture some bacterial patterns which can be reproduced by the model \eqref{(4)}. More precisely, we focus on bacterial honeycomb pattern, reported for example in \cite{Thar} , and symmetrical spot pattern  generated for instance by \textit{Escherichia coli}  \cite{budrene2}.

In this subsection, $D_u=0.0625$,  $\lambda=32$ and $f(u)=u^2(1-u)$. Initial conditions and computational domain  are the same as in Subsection 4.2. The scheme adopted  is similar to \eqref{x1}--\eqref{x2}, but since the logistic source term is now cubic, $\widetilde u_K^{n+1}$ in \eqref{x1}  will be replaced by $\vert\widetilde u_K^{n+1}\vert ^2$. The time-step size used for computations is $\Delta t=10^{-1}$.

\begin{figure}[h!]
\subfigure{\includegraphics[width=6.cm]{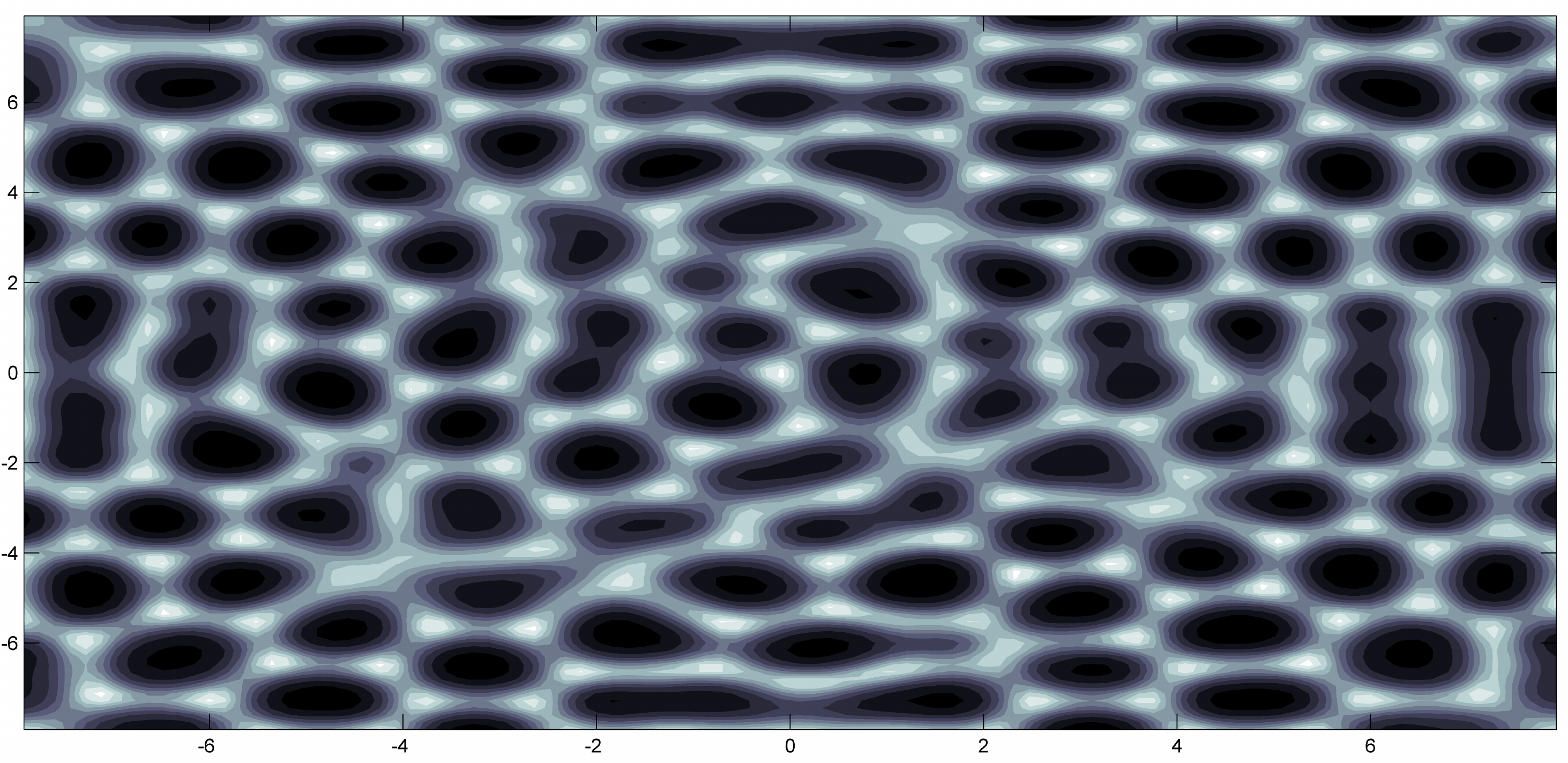}}
\subfigure{\includegraphics[width=6.2 cm]{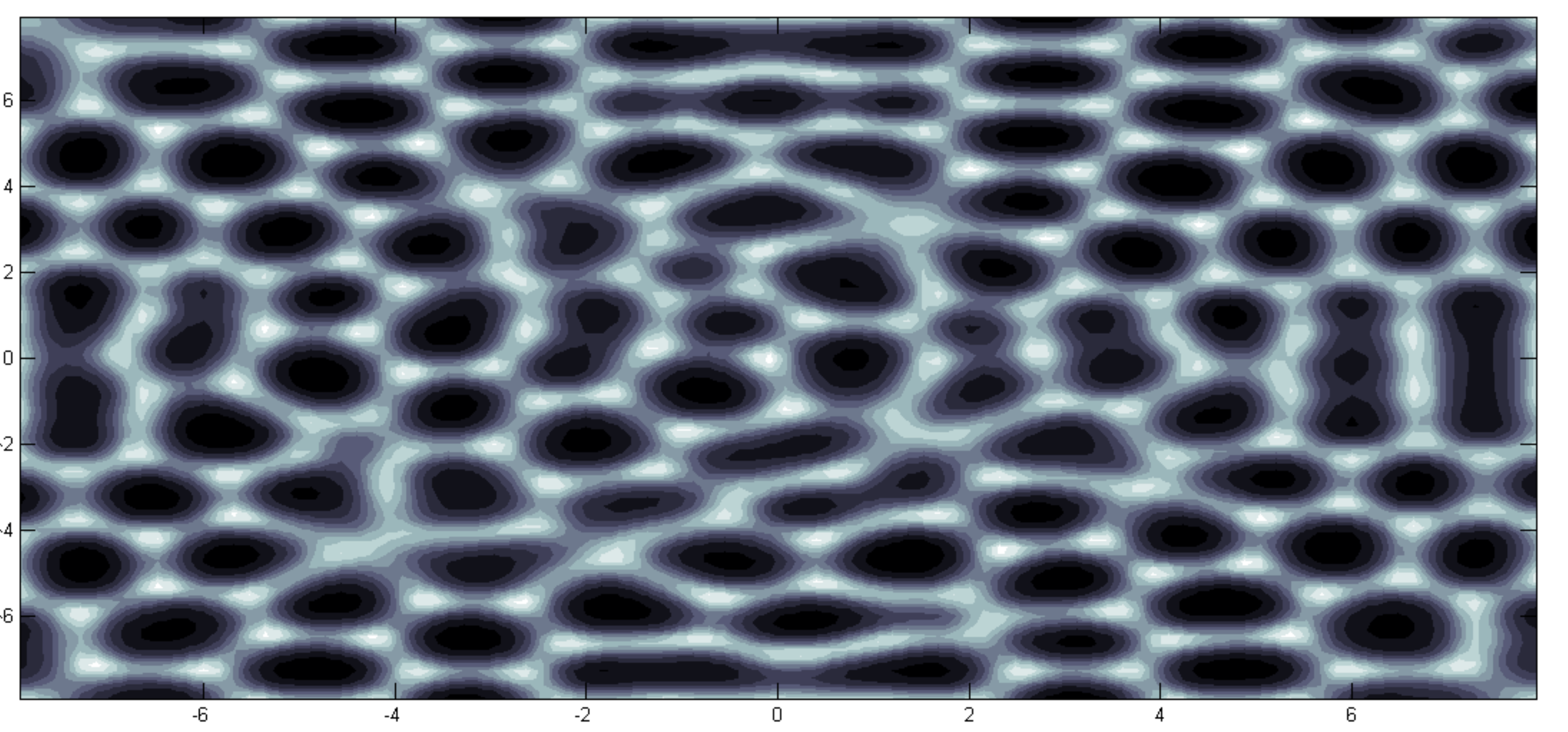}}
\caption{Solution $(u)$ of  \eqref{(4)} with cubic growth at $t=150$  computed via the SSTLI method (left) and the ESTLI method (right) with $\chi=6$ (black indicates low cell density, white indicates high cell density) }
\label{Fig.11}
\end{figure}

\begin{figure}[h!]
\subfigure{\includegraphics[width=6.cm]{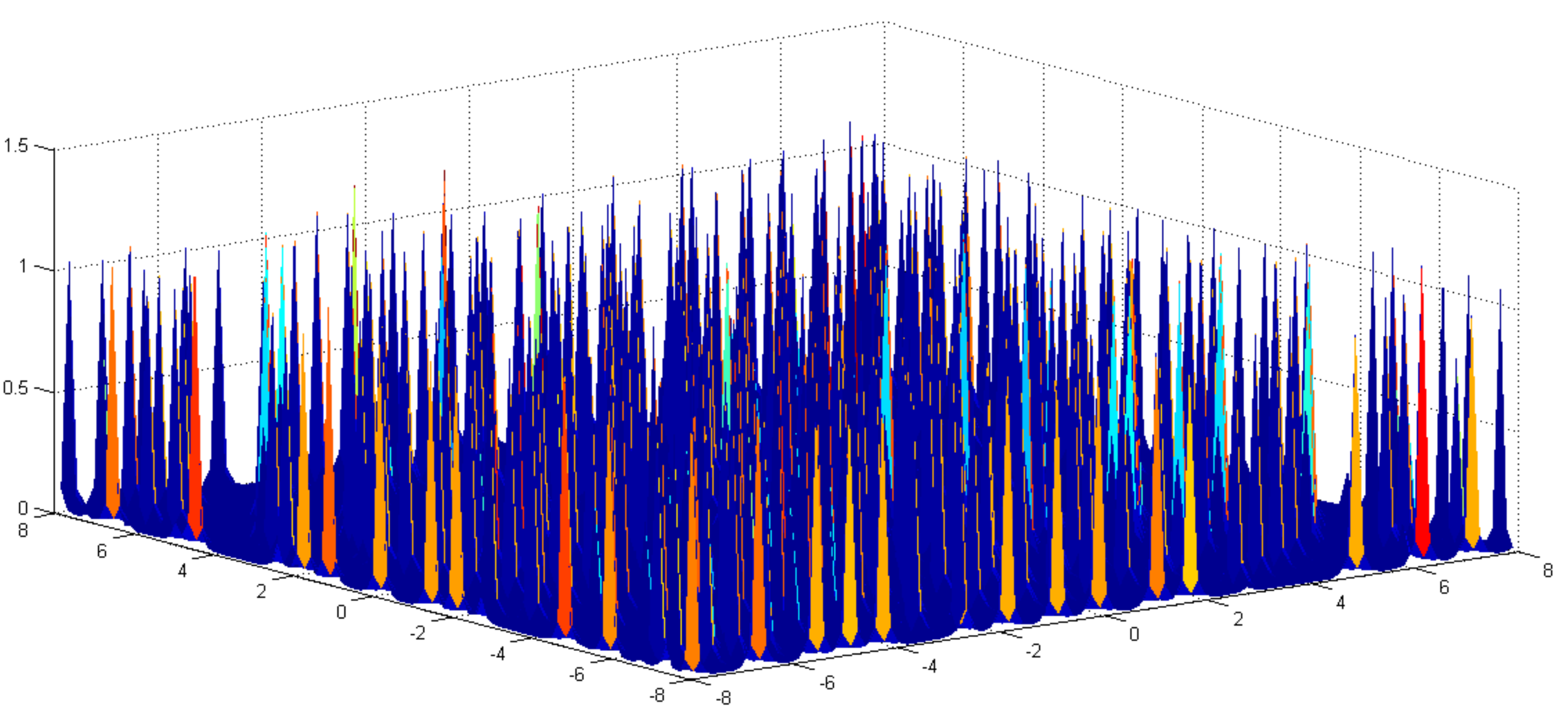}}
\subfigure{\includegraphics[width=6.cm]{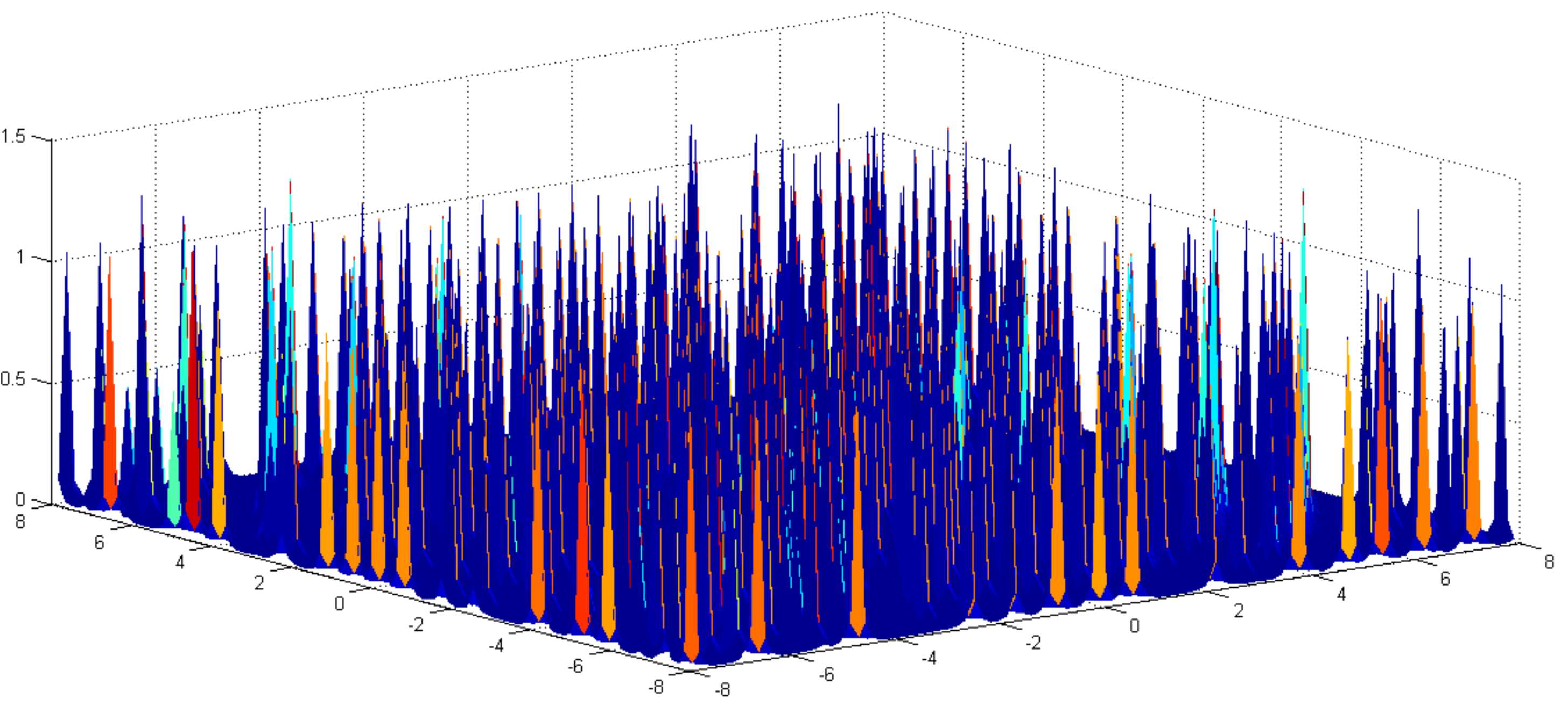}}
\caption{Three-dimensional plot of the solution $(u)$ of  \eqref{(4)} with cubic growth at $t=10$  computed via the SSTLI method (left) and the ESTLI method (right) with $\chi=120$ }
\label{Fig.12}
\end{figure}

\begin{figure}[h!]
\subfigure{\includegraphics[width=6.cm]{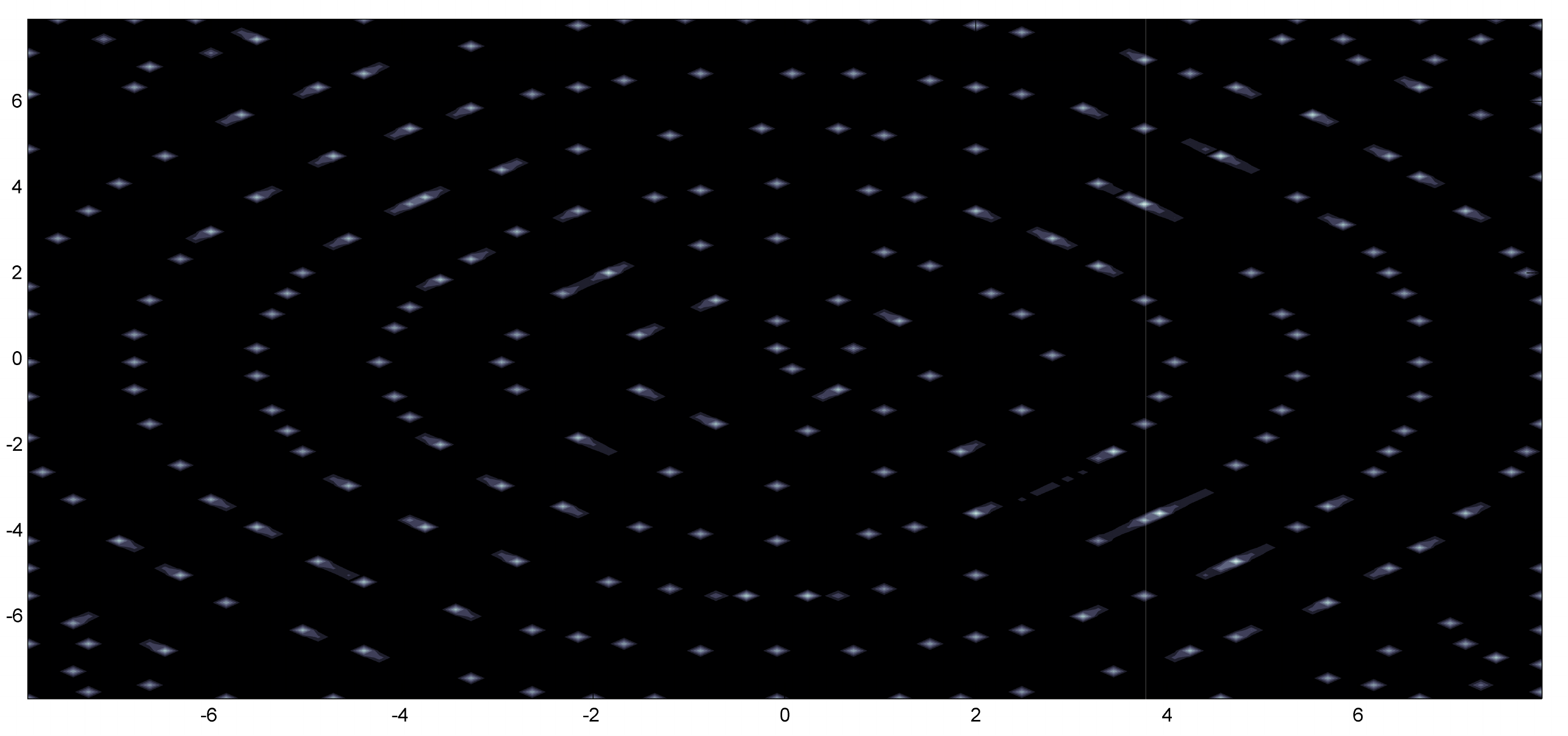}}
\subfigure{\includegraphics[width=6.cm]{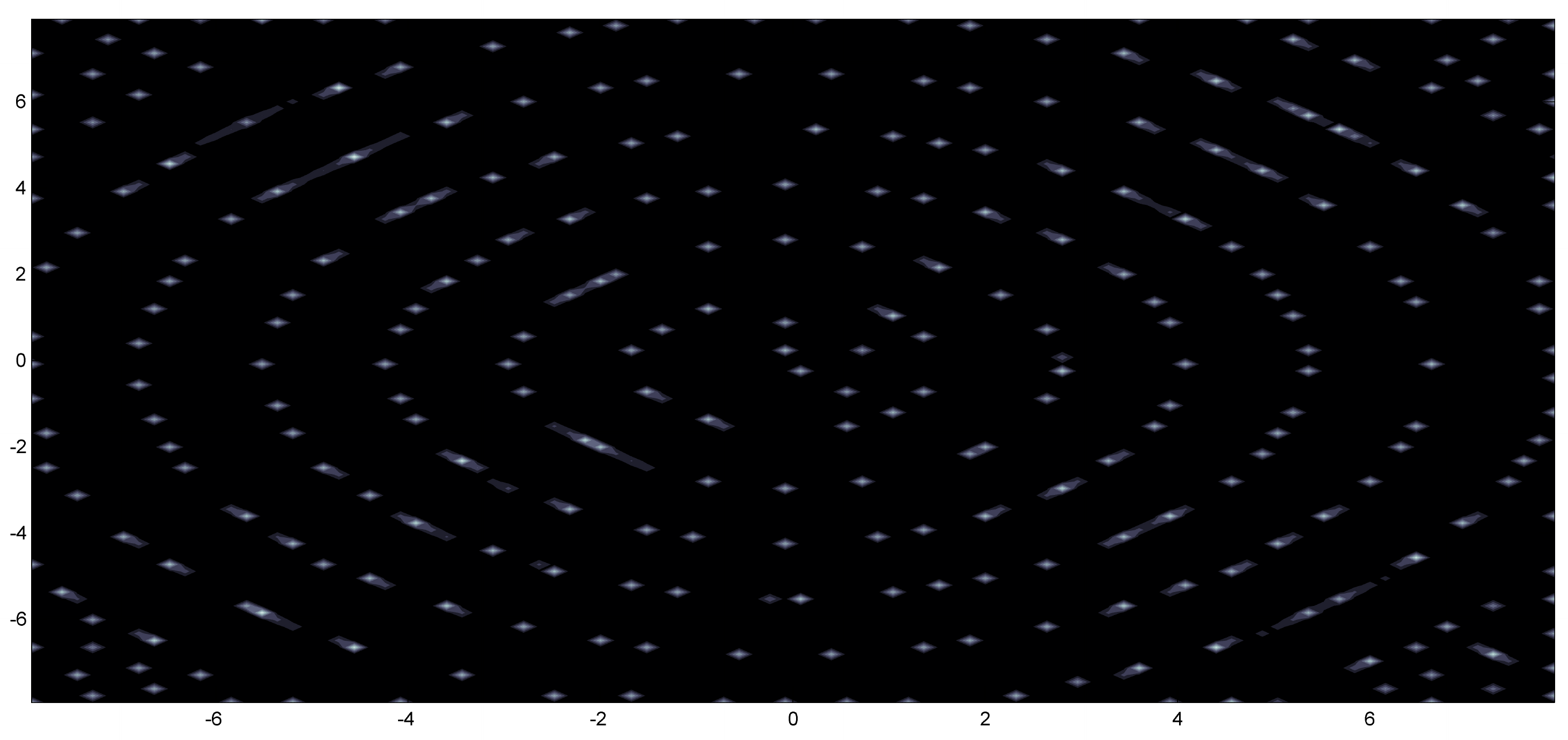}}
\caption{Solution $(u)$ of  \eqref{(4)} with cubic growth at $t=10$  computed via the SSTLI method (left) and the ESTLI method (right) with $\chi=120$ (black indicates low cell density, white indicates high cell density) }
\label{Fig.13}
\end{figure}

For $\chi=6$, the computed solutions at $t=150$ are shown in Fig. \ref{Fig.11}. The honeycomb pattern is observed for both ESTLI and SSTLI numerical solutions. When $\chi=120$, we see from Fig. \ref{Fig.12} that the computed solutions exhibit a very spiky behavior at $t=10$, which means that spot pattern appears (see Fig. \ref{Fig.13}). Despite this fact, the two methods remain positivity preserving which confirms their robustness and reliability.

\section{Conclusion}

In this paper, two linearized implicit methods relying on a single-layer neural network are developed to solve  two-dimensional Keller-Segel systems: ESTLI and SSTLI methods. In terms of accuracy, the numerical tests performed demonstrate the superiority of the developed methods on semi-implicit method. However, concerning computational cost, only SSTLI method is nearly as efficient as semi-implicit scheme for time-consuming numerical tests. Moreover, the two numerical methods are robust and easily applicable to any problem which can be solved by the semi-implicit approach.

\end{document}